# NUMERICAL APPROXIMATION OF LAMBERT W FUNCTION FOR REAL VALUES BY UNIQUE METHOD OF QUADRATIC APPROXIMATION

**Narinder Kumar Wadhawan,**

Civil Servant, Indian Administrative Service Now Retired, House No. 563, Sector 2, Panchkula-134109, Haryana, India, E-mail: narinderkw@gmail.com

*Abstract*: This paper introduces a novel method for the numerical approximation of the Lambert W function, $W(x) = y$, in the real domain. By linearly approximating the natural logarithm, the function is transformed into a quadratic equation whose roots refine the initial approximation. Iterative solving of this equation yields high precision, with iterations determined by the desired accuracy. Two methods for positive $x$ are proposed, the first expresses the function as $(z_1 + a_1) \ln(z_1 + a_1) = x$, iteratively refining the approximation of $z$, where $z_1$, is the initial approximation and $a_1 \ll z_1$, while the second method is based on $\ln(y_1 + a_1) + y_1 + a_1 = \ln(x)$, where $y_1$ is the initial approximation, and $a_1 \ll y_1$. For $x$ between $0$ and $1/e$, the method is extended to approximate $W(-x) = -y$. Unlike Newton or Halley's methods, this approach handles both branches without constraints on initial assumptions, covering a broad range of values. Extensive examples and a software algorithm validate the method's accuracy, offering a precise and flexible tool for numerical analysis.

## 1. Introduction
The Lambert W function, denoted as $W(x)$, is the inverse function of $y = x e^x$. It was first explored by Johann Heinrich Lambert in 1758, who provided a series solution to the equation $xe^x = y$, which could be written
$$W(x) = y. \qquad (1.1)$$
He later refined it by including higher-order terms [17], [18]. In 1784, Leonhard Euler extended Lambert's work by deriving a series solution for a similar equation [12]. Euler simplified the equation by taking the limit as $x$ tended to infinity, ultimately arriving at a convergent series for $W(x)$.

The Lambert W function is multivalued, meaning it has multiple branches. For real numbers, the function is solvable when $x > -1/e$, with the principal branch $W_0(x)$ valid for $x > -1/e$ and another branch $W_{-1}(x)$ existing for $-1/e \leq x < 0$. These two branches are of particular importance in many scientific applications.

In recent years, the Lambert W function has found applications across diverse fields such as fluid dynamics, where it is used to model flow in porous media, and chemical engineering [5], [7], in which it helps solve reaction rate equations. In neuroimaging [27], it aids in analysing nonlinear signal transformations, and in epidemiology [30], it is employed in modelling the spread of infectious diseases. A recent study in material science [6] applied the Lambert W function to describe phase transitions in complex polymers. The function also plays a crucial role in quantum mechanics, particularly in finding exact solutions to the Schrödinger equation [15] and expression of the corrections to the dispersion relations of giant magnons, single spikes and GKP strings [13]. It can solve equation linked to Bernoulli numbers and Todd genus by means of the two real branches $W_0$ and $W_{-1}$ [16]. The Lambert W function is widely used to solve equations of the form $xe^x = y$ and other equations that can be transformed into this form. In addition to the application already explained, use of the Lambert W function spans various disciplines, including, porous media [10], statistics, [23], Michaelis-Menten kinetics [21], time-dependent flow in simple branch hydraulic systems [24], crystal growth [2], pooling of tests for infectious diseases [3] [8] [14], exact solution of QCD coupling constant [11], resonance of the delta-shell potential [19], phase separation of polymer mixtures [4], electromagnetic surface wave propagation [20], and viscous flow [21]. It also finds use in thermodynamic equilibrium, Wein's displacement law [9] and Quantum-mechanical double-well Dirac delta function model for equal charges [26].

In this paper, a novel method for the numerical approximation of both branches of the Lambert W function is presented. Unlike traditional methods such as Newton's method [1], [29] or Halley's method [25], the approach used here is based on quadratic approximation. Numerical approximations of both the branches of the function have been made by using a unique method that transforms a real number say $y$ into exponential form and vice versa as has been proved in a published paper [28]. Based upon the theory, following identities were proved [28],

$$y/(y-1) = e^{2/(2y-1)}, \tag{1.2}$$

$$\ln y = 2 \sum_{n=2}^{y} 1/(2y-1), \tag{1.3}$$

*provided $y$ is appreciably large.* These identities will be used to evaluate the values of Lambert W function. Since my assumption requires $y$ to be large and it is not always possible to have $y$ with large value, therefore, repeated corrections using recursive relation,

$$y_{n+1} = y_n + a_n, \tag{1.4}$$

will be applied to achieve precise approximations of $y$ that will finally equal $y_{n+1}$, where $a_n$ is the root of the nth resultant quadratic equation in variable $a_n$. In this paper, corrections three to four times, have been applied but more corrections may be needed to enhance precision based on accuracy requirements. Unlike Newton's and Halley's methods, this quadratic equation approximation does not require the differentiability or continuity of the function, nor does it impose conditions on the function's derivative at assigned values. This versatility allows the method to handle incorrect initial assumptions without failing, making it applicable to a broader range of problems without the need for preprocessing steps to refine guesses—issues that plague existing approximation methods.

The method features a self-correcting capability, minimising human error, whereas traditional methods can be hindered by initial assumption errors, wasting time and resources. A significant advantage is its ability to converge within three to four iterations, regardless of the initial guess, compared to other methods that require multiple iterations for high accuracy. Further details are given in section 8 where a companion has been made. The validation of this method includes solving both branches of the Lambert W function, with comprehensive examples presented in tabular form. Additionally, an algorithm along with its pseudocode are provided for the numerical approximation of $W(x)$ for specified values, serving as a practical computational tool.

## 2. Approximation, Convergence and Rate of Convergence

In this section, theory of approximation, convergence of iterated results and associated error will be taken up. Results derived in this section will be verified in the sections that follow.

### 2.1 Quadratic Approximation of Lambert W Function $W(x) = y$

Assuming $y$ equal to $y_1 + a_1$, where $y_1$ is assigned an arbitrary initial value, the equation $y\, e^y = x$ takes the form $y_1 + a_1 + \ln(y_1 + a_1) = \ln(x)$. This equation is then transformed to a quadratic equation in $a_1$ and the root of this equation yields values $a_1$ and the value of $y$ is approximated as $y_1 + a_1$. In some cases, where $y$ is considered equal to $\ln(z)$, the equation $y\, e^y = x$, then takes the form $z \ln(z) = x$ and $z$ is assumed equal to $z_1 + a_1$, where $z_1$ is assigned an arbitrary initial value. The equation $(z_1 + a_1) \ln(z_1 + a_1) = x$ is then transformed to a quadratic equation in $a_1$ and its roots approximate the value of $z = z_1 + a_1$. With this background, the detailed method is explained hereinafter.

It is proved in a published paper [28] that the exponential term $e^{2/(2y-1-1/y^3)}$ can be approximated with $y/(y-1)$ provided value of $y$ is large, otherwise it will yield rough approximation. If $y$ is assumed as $-y$, term $e^{-2/(2y+1-1/y^3)}$ will approximate term $y/(y+1)$. Taking logarithm, term $2/(2y+1-1/y^3)$ will approximate term $ln\{(y+1)/y\}$.

*Before continuing further, I clarify the use of signs of approximation ($\simeq$) and equality ($=$). Lambert W Function involves exponential terms which are irrational and that needs use of approximation sign in stead of equality but as the paper proceeds, it will be observed, equality sign has only been used. Such usages do not imply that results are exact but on the other hand, mean, these are approximations of the actual results.*

Coming to the approximation of $2/(2y+1-1/y^3)$ to $ln\{(y+1)/y\}$, when $y$ is a large quantity, term $1/y^3$ can be ignored being small, then

$$ln\{(y+1)/y\} \simeq 2/(2y+1). \tag{2.1}$$

This Equation (2.1) will invariably be used in this paper. Lambert W function $W(x) = y$, can also be written $y\, e^y = x$ or $z\, ln(z) = x$, where $y = ln(z)$ and $z$ is a real number neither equal to zero nor infinity. Assuming $z = z_1 + a_1$, it can be written

$$ln\, z = ln(z_1 + a_1) = ln\, z_1 + ln\left(\frac{z_1/a_1 + 1}{z_1/a_1}\right). \tag{2.2}$$

Simplifying Equation (2.2) with the use of Equation (2.1),

$$ln\, z \simeq ln\, z_1 + 2a_1/(a_1 + 2z_1). \tag{2.3}$$

Therefore,

$$(a_1 + z_1)\, ln\, z_1 + 2a_1(a_1 + z_1)/(a_1 + 2z_1) \simeq x. \tag{2.4}$$

Multiplying both sides by $(a_1 + 2z_1)$ and simplifying, a quadratic equation in $a_1$ is obtained. This equation can be written as

$$a_1^2 - l_1 a_1 - m_1 \simeq 0, \tag{2.5}$$

where coefficients,

$$l_1 = -(3z_1\, ln\, z_1 + 2z_1 - x)/(ln\, z_1 + 2),$$
$$m_1 = 2z_1(x - z_1\, ln\, z_1)/(ln\, z_1 + 2),$$

Root $a_1 = l_1/2 + 1/2\sqrt{l_1^2 + 4m_1}$ of the quadratic equation will be used, *since that finally gives positive real value of $z_n + a_n$, where $z_{n+1} = z_n + a_n$*. It is essential to have a positive value of $z_n + a_n$ otherwise, negative value will lead to natural logarithm of a negative quantity and logarithm of a negative quantity does not exist in real domain.

Once the roots are determined, the positive value of $z_1 + a_1$ say $z_2$ will approximate $z$. For obtaining fine approximation of $z$, $z_3$ will then be considered equal to $z_2 + a_2$ and value of $a_2$ will be determined following the same procedure as that of $a_1$. When $z_2 + a_2$ has been determined, the process of approximation will be continued till $z_n + a_n$ is determined, where $n$ is number of times the correction is applied to obtain precise approximation. It should also be noted that $z_1$ should be assigned a value so that $a_n$ should have a real value after iteration and $z_n + a_n$ should also have a real positive value after iteration, but it should not be zero or infinity. How such a choice should be made, has been explained and illustrated with exhaustive examples in forthcoming sections.

In some cases, approximation of value $y$ has been made by taking natural logarithm of both sides of Lambert W Function, $ye^y = x$. In such cases, $y + ln\, y = ln\, x$. Thereafter, $y$ is considered equal to $y_1 + a_1$ and the equation takes the form $(y_1 + a_1) + ln(y_1 + a_1) = ln\, x$. This equation can be written as $(y_1 + a_1) + ln\, y_1 + ln(1 + y_1/a_1) = ln\, x$. Applying Equation (2.1) derived in a published paper [28], $ln(1 + y_1/a_1)$ can be written $2a_1/(a_1 + 2y_1)$ and equation, $y + ln\, y = ln\, x$ takes the form .

$$(y_1 + a_1) + ln\, y_1 + 2a_1/(a_1 + 2y_1) = ln\, x. \tag{2.6}$$

This equation on simplification is written in quadratic form in variable $a_1$,
$$a_1^2 - l_1 a - m_1 = 0,$$
which has its coefficients,
$$l_1 = -\{3y_1 + 2 - \ln(x/y_1)\},$$
$$m_1 = -2y_1\{y_1 - \ln(x/y_1)\}.$$
Out of the two roots of the quadratic equation, the root $a_1 = (1/2)\left(l_1 + \sqrt{l_1^2 + 4m_1}\right)$ will be used. For determining precise approximation, correction is applied $n$ times to obtain $y_n + a_n$, which will be equal to precise value of $y$.

Above said, quadratic approximations will be used to evaluate $y$, when value of $x$ is given in equation $ye^y = x$.

### 2.2. Convergence First Method
#### 2.2.1. Convergence of $\ln z = y$

In this section, convergence of $\ln z = y$ will first be proved and thereafter, convergence of Lambert W Function, $z \ln z = x$, will be taken up. In the first instance, considering $z = z_1 + a_1$ yields $\ln(z_1 + a_1) = y$, where $z$ is a positive real number and $z_1$ is an arbitrarily assigned positive real number. This equation on simplification yields $\ln z \simeq \ln z_1 + 2a_1/(a_1 + 2z_1) = y$ [28] and value of $a_1$ is thus given by equation

$$a_1 = 2z_1(y - \ln z_1)/(2 - y + \ln z_1). \tag{2.7}$$

Considering value of $z_2 = z_1 + a_1$ and on simplification, yields

$$z_2/z_1 = -1 + 4/(2 - y + \ln z_1). \tag{2.8}$$

Considering value of $z_3 = z_2 + a_2$ and on simplification, yields

$$z_3/z_2 = -1 + 4/(2 - y + \ln z_2), \tag{2.9}$$

and finally, value of $z_{n+1} = z_n + a_n$ is given by relation,

$$z_{n+1}/z_n = -1 + 4/(2 - y + \ln z_n), \tag{2.10}$$

where $n$ denotes the number of times the correction is applied. In ideal case, $n$ is that number which yields $z_{n+1}/z_n = 1$. For approximation of $\ln z$, value to $z_1$ is arbitrarily assigned so that values of $z_2, z_3, z_4, \ldots, z_{n+1}$ obtained are positive and real. What value should be assigned to $z_1$ has been explained in sections 2.2.1.1 and 2.2.1.2. Assignment of value to $z_1$ will lead to two types of cases.

*2.2.1.1. First Type, when $\ln z_1 - y$ is positive but less than 2:* Referring to Equation (2.8), denominator $2 - y + \ln z_1$ will have a value more than 2, when $z_1$ is assigned a value such that $\ln z_1 > y$. That results in $z_2 < z_1$ from Equation (2.8). *Since Equation (2.8) has been satisfied, $\ln(z_1 + a_1)$ will aapriximate $y$ more precisely than $\ln z_1$.*

Referring to Equation (2.9) and $z_2 < z_1$, denominator $2 - y + \ln z_2$ will either be more than 2 or equal to 2 or even less than 2. If it is equal to 2, no further correction is required and the result will be $\ln z_2 = y$. It is submitted such a situation will only arise, when assumption of $z_1$ is such that $\ln z_1$ closely approximates $y$ otherwise magnitude of $2 - y + \ln z_2$ will be more than 2 but less than $2 - y + \ln z_1$. Since Equation (2.9) has been satisfied, therefore, $|z_3 - z_2|$ will be less than $|z_2 - z_1|$ and $\ln(z_2 + a_2)$ *will aapriximate $y$ more precisely than $\ln(z_1 + a_1)$.* In this way, with successive iterations, $|z_4 - z_3| < |z_3 - z_2|$, $|z_5 - z_4| < |z_4 - z_3|, \ldots, |z_{n+1} - z_n| < |z_n - z_{n-1}|$.

*2.2.1.2. Second Type, when $\ln z_1 - y$ is negative but $2 + \ln z_1 - y$ is positive quantity:* In this type of cases, $(2 - y + \ln z_1)$ is less than 2 but more than 0 and satisfaction of Equation (2.8), yields $z_2 - z_1$ a positive quantity or $z_2 > z_1$ and $\ln(z_1 + a_1)$ *will aapriximate $y$ more precisely than $\ln z_1$.*

Referring to Equation (2.9) and $z_2 > z_1$, denominator $2 - y + \ln z_2$ will either be less than 2 or equal to 2 or even more than 2. If it is equal to 2, no further correction is required and the result will be $\ln z_2 = y$. Such a situation will only arise, when assumption of $z_1$ is such that

$\ln z_1$ closely approximates $y$ otherwise magnitude of $2 - y + \ln z_2$ will be less than 2 but more than $2 - y + \ln z_1$. Since Equation (2.9) has been satisfied, therefore, magnitude of $z_3 - z_2$ written as $|z_3 - z_2|$ will be less than $|z_2 - z_1|$ and $\ln(z_2 + a_2)$ *will aapriximate $y$ more precisely than* $\ln(z_1 + a_1)$. In this way, with successive iterations, $|z_4 - z_3| < |z_3 - z_2|$, $|z_5 - z_4| < |z_4 - z_3|, \ldots, |z_{n+1} - z_n| < |z_n - z_{n-1}|$. Difference in this type of cases from first type, is that here $\ln z_1 - y$ which was initially negative, with each iteration get less negative and finally will reach zero.

Considering magnitudes of $z_2 - z_1, z_3 - z_2, z_4 - z_3, \ldots, z_{n+1} - z_n$ and the position explained above, it proves

$$|z_2 - z_1| > |z_3 - z_2| > |z_4 - z_3| > \cdots > |z_{n+1} - z_n|, \qquad (2.11)$$

or

$$|a_1| > |a_2| > |a_3| > \cdots > |a_n|,$$

or $a_1 + a_2 + a_3 + \cdots + a_n$ is convergent and when $n$ is appreciably large, $z_{n+1}/z_n$ tends to unity. It is submitted, $z_1$ is not to be assigned a value that results in $-1 + 4/(2 - y + \ln z_1)$ zero or a negative quantity as logarithm of zero or a negative quantity is not determinable in real domain. In other words, $-2 < -y + \ln z_1 < 2$.

*Example:* Let $y = 0.8$, in equation $\ln z = y$. Assuming $z_1 = 3$, and applying equations (2.7) to (2.11), values of $z_2, z_3, z_4, \ldots, z_{n+1}$ obtained are

$z_2 = 2.220541132, \quad z_3 = 2.225540931, \quad z_4 = 2.225540928, \quad z_5 = 2.225540928.$

For the result correct up to nine places of decimals, it is observed, $z_5/z_4$ precisely approximates one, therefore, $z_5$ is the approximated value that yields $\ln z_5 = 0.8$. In terms of differences,

$|z_2 - z_1| = 0.779458868, \qquad |z_3 - z_2| = 0.004999799,$
$|z_4 - z_3| = 3 \times 10^{-9}, \qquad |z_5 - z_4| = 0.000000000$

and these are in decreasing order. Therefore,

$$a_1 + a_2 + a_3 + \cdots + a_n,$$

is a convergent series and $z_1 + a_1 + a_2 + a_3 + \cdots + a_n$ sums up to $z_{n+1}$, where $n$ is 5.

### 2.2.2. Convergence of Lambert W Function using quadratic approximation

It has been proved that positive value of $z$ in equation $\ln z = y$ obtained by approximation converges to $z_1 + a_1 + a_2 + a_3 + \cdots + a_n$, where $z_1$ is arbitrarily assumed positive real number and

$$|a_1| > |a_2| > |a_3| > \cdots |a_{n-1}| > |a_n|.$$

Based on above proof, following deductions are made.

1) Since $z$ in equation, $\ln(z) = y$, converges to a positive nonzero value $z_{n+1}$, therefore, $z$ in product of $z$ and $\ln z$ equal to $x$ or in $z \ln z = x$, will also coverage to a positive value.
2) Since $z$ in equation, $\ln(z) = y$, converges to a positive nonzero value $z_{n+1}$, therefore, $z$ in division of $\ln z$ by $z$ equal to $x$ or in $(1/z) \ln z = x$, will also coverage to a positive value.
3) Since $z$ in equation, $\ln(z) = y$, converges to a positive nonzero value $z_{n+1}$, therefore, $z$ in sum of $\ln z$ and $z$ equal to $x$ or in $z + \ln z = x$, will also coverage to a positive value.
4) Since $z$ in equation, $\ln(z) = y$, converges to a positive nonzero value $z_{n+1}$, therefore, $z$ in algebraic sum of $\ln z$ and $-z$ equal to $x$ or in $-z + \ln z = x$, will also coverage to a positive value.

### 2.3. Convergence Second Method

#### 2.3.1. Convergence of Lambert W Function Expressed As $z \ln(z) = x$, Where $z = e^y$

Referring to section 3.1, Lambert W Function is written as $z \ln z = x$, where $z = e^y$.

Considering $z = z_1 + a_1$, it gets transformed to a quadratic equation $a_1^2 - l_1 a - m_1 = 0$, where
$$l_1 = -(3z_1 \ln z_1 + 2z_1 - x)/(\ln z_1 + 2),$$
$$m_1 = 2z_1(x - z_1 \ln z_1)/(\ln z_1 + 2),$$
and relevant root of quadratic equation is
$$a_1 = (1/2)\left(l_1 + \sqrt{l_1^2 + 4m_1}\right).$$
Then
$$z_2 = z_1 + a_1 = (1/2)\left(l_1 + \sqrt{l_1^2 + 4m_1}\right).$$
For precise approximation, it is further iterated and value of $z_3 = (z_2 + a_2)$, $z_4 = (z_3 + a_3), \ldots, z_{n+1} = (z_n + a_n)$ are found. Detail procedure is mentioned in section 3.1 and is not repeated here. It is reiterated, quadratic approximation is based on formula
$$\ln (y + 1)/y\} \simeq 2/(2y + 1),$$
which has already been proved in a published paper [28], therefore, value of $z_2 = z_1 + a_1$ will satisfy more precisely equation $z \ln z = x$ as compared to $z = z_1$. Similarly, $z_3$ will satisfy more precisely equation $z \ln z = x$ than $z = z_2$ and finally $z_{n+1} = z_n + a_n$ will satisfy more precisely equation $z \ln z = x$ than $z_n$. That leads to result
$$|z_1 \ln z_1 - x| > |z_2 \ln z_2 - x| > |z_2 \ln z_2 - x| > \cdots > |z_n \ln z_n - x|. \quad (2.12)$$
Since according to Equation (3.7), $l_2 = -(3z_2 \ln z_2 + 2z_2 - x)/(\ln z_2 + 2)$, and $z_2 \ln z_2$ approximates $x$, considering $z_1 \geq 1$, then $l_2$ is a negative quantity and $m_2 = 2z_2(x - z_2 \ln z_1)/(\ln z_2 + 2)$ can have positive or negative value depending upon assumption of value of $z_1$. That is whether $z_2 \ln z_2$ is less or more than $x$ but since $z_2 \ln z_2$ approximates $x$, therefore, quantity $|x - z_2 \ln z_2|$ will be small as compared to $|3z_2 \ln z_2 + 2z_2 - x|$. Let $l_2 = -L_2$, then on simplification,
$$a_2 = (L_2/2)\{-1 + (1 + 4m_2/L_2^2)^{1/2}\} = m_2/L_2. \quad (2.13)$$
Since $z_2 \ln z_2$ approximates $x$, therefore, $a_2$ on substitution of the values of $L_2$, $m_2$ and simplification, is given by equation,
$$a_2 = 2z_2(x - z_2 \ln z_1)/(3z_2 \ln z_2 + 2z_2 - x) = (x - z_2 \ln z_2)/(\ln z_2 + 1). \quad (2.14)$$
Similarly,
$$a_3 = (x - z_3 \ln z_3)/(\ln z_3 + 1), \quad (2.15)$$
$$\ldots$$
$$a_n = (x - z_n \ln z_n)/(\ln z_n + 1). \quad (2.16)$$
According to inequality (2.12) and in view of the fact, slight change in values of $z_2, z_3, z_4, \ldots, z_{n+1}$ on iteration does not affect appreciably the values of denominators $(1 + \ln z_2)$, $(1 + \ln z_3)$, $(1 + \ln z_4), \ldots, (1 + \ln z_n)$, therefore,
$$|a_2| > |a_3| > |a_4| > \cdots > |a_n|,$$
and thus, value of $z$ converges to $z_{n+1}$, where
$$z_{n+1} = z_1 + a_1 + a_2 + a_3 + \cdots + a_n,$$
value of $z_1$ is arbitrarily assumed, value of $a_n = (1/2)\left(l_n + \sqrt{l_n^2 + 4m_n}\right)$ and value of $l_n$ and $m_n$ are given by Equations (3.10) and (3.11).

2.3.2. Convergence of Lambert W Function Expressed as $(1/z) \ln (z) = x$, Where $z = e^y$

Following the same procedure as explained in section 2.3.1, it can be proved,
$$|a_2| > |a_3| > |a_4| > \cdots > |a_n|,$$
and thus, value of $z$ converges to $z_{n+1}$, where
$$z_{n+1} = z_1 + a_1 + a_2 + a_3 + \cdots + a_n,$$
value of $z_1$ is arbitrarily assumed, value of $a_n = (1/2)\left(l_n + \sqrt{l_n^2 + 4m_n}\right)$ and value of $l_n$ and $m_n$ are given by equations,
$$l_n = -(3z_n x - \ln z_n - 2)/x,$$
$$m_n = 2z_n (\ln z_n - z_n x)/x.$$

Since Lambert W Function in this case, pertains to negative branch, therefore, it can be proved, its other value say $z'$ also converges to
$$z'_{n+1} = z_1 + a'_1 + a'_2 + a'_3 + \cdots + a'_n,$$
where value of $a'_n = (1/2)\left(l'_n - \sqrt{l'^2_n + 4m'_n}\right)$, $|a'_2| > |a'_3| > |a'_4| > \cdots > |a'_n|$ and value of $l'_n$ and $m'_n$ are given by equations.
$$l'_n = -(3z'_n x - \ln z'_n - 2)/x,$$
$$m'_n = 2z'_n (\ln z'_n - z'_n x)/x.$$

### 2.3.3. Convergence of Lambert W Function Expressed As $\ln(y) + y = \ln(x)$

Referring to section 3.2, Lambert W Function $ye^y = x$ is written as $\ln(y) + y = \ln(x)$ and continuing as explained in section 2.3.1, value of $a_2$ is given by equation $a_2 = m_2/L_2$, where $L_2 = -l_2$ On substituting the values and simplifying it,
$$a_2 = -2y_2\{y_2 - \ln(x/y_2)\}/\{3y_2 + 2 - \ln(x/y_2)\} = \{\ln(x/y_2) - y_2\}/(1 + 1/y_2) \quad (2.17)$$
Similarly,
$$a_3 = \{\ln(x/y_3) - y_3\}/(1 + 1/y_3), \quad (2.18)$$
$$\ldots$$
$$a_n = \{\ln(x/y_n) - y_n\}/(1 + 1/y_n). \quad (2.19)$$
Since $y_2$ approximates more precisely the equation $\ln(y) + y = \ln(x)$ than $y_1$ and $y_3$ approximates still more precisely the equation $\ln(y) + y = \ln(x)$ than $y_2$ and so on. This leads to result
$$|\ln(x/y_1) - y_1| > |\ln(x/y_2) - y_2| > |\ln(x/y_3) - y_3| > \cdots > |\ln(x/y_n) - y_n| \quad (2.20)$$
According to inequality (2.20) and in view of the fact, there is extraordinarily slight changes in the values of $y_2, y_3, y_4, \ldots, y_{n+1}$ and that does not affect appreciably the values of denominators $(1 + 1/y_2)$, $(1 + 1/y_3), (1 + 1/y_4), \ldots, (1 + 1/y_n)$, therefore,
$$|a_2| > |a_3| > |a_4| > \cdots > |a_n|,$$
and thus, value of $y$ converges to $y_{n+1}$ where
$$y_{n+1} = y_1 + a_1 + a_2 + a_3 + \cdots + a_n,$$
value of $y_1$ is arbitrarily assigned, value of $a_n = (1/2)\left(l_n + \sqrt{l^2_n + 4m_n}\right)$ and value of $l_n$ and $m_n$ are given by Equations (3.22) and (3.23).

### 2.3.4. Convergence of Lambert W Function Expressed As $\ln(y) - y = \ln(x)$

Following the same procedure as explained in section 2.3.2, it can be proved,
$$|a_2| > |a_3| > |a_4| > \cdots > |a_n|,$$
and thus, value of $y$ converges to $y_{n+1}$ where
$$y_{n+1} = y_1 + a_1 + a_2 + a_3 + \cdots + a_n,$$
value of $y_1$ is arbitrarily assigned, value of $a_n = (1/2)\left(l_n + \sqrt{l^2_n + 4m_n}\right)$ and values of $l_1$ and $m_1$ are given by equations
$$l_n = -\{3y_n - 2 + \ln(x/y_n)\},$$
$$m_n = -2y_n\{y_n + \ln(x/y_n)\}.$$
Since Lambert W Function in this case, pertains to negative branch, therefore, it can be proved, its other value say $y'$ also converges to
$$y'_{n+1} = y_1 + a'_1 + a'_2 + a'_3 + \cdots + a'_n,$$
where $|a'_2| > |a'_3| > |a'_4| > \cdots > |a'_n|$, value of $a'_n = (1/2)\left(l'_n - \sqrt{l'^2_n + 4m'_n}\right)$ and value of $l'_n$ and $m'_n$ are given by equations,
$$l'_n = -\{3y'_n - 2 + \ln(x/y'_n)\},$$
$$m'_n = -2y'_n\{y'_n + \ln(x/y'_n)\}.$$

### 2,3.5. Graphical Presentation of Convergence

To further corroborate convergence of quadratic approximation, I refer to the data given in Table 3.1. Figures 1, 2, 3, 4 and 5 are drawn according to the data given in rows at serial numbers 2, 3, 4, 5 and 6 respectively in Table 3.1. In these Figures, along x-axis, numerals 1, 2, 3, 4, 5 and 6 correspond to value of $z_1, z_2, z_3, z_4, z_5$ and $z_6$, where $z_1$ is an assumed value. Applying correction four times results in approximation precise to ten decimal points.

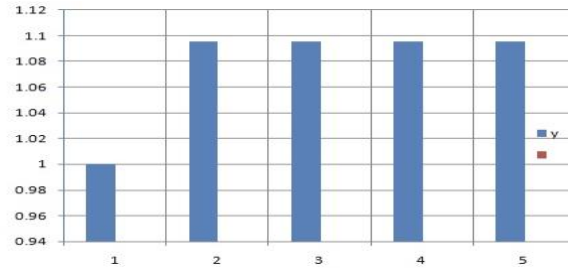

Figure 1 Convergence of $z_1 + a_1 + a_2 + a_3 + a_4$, when $x = 0.1$ and $z_1$ assumed is 1

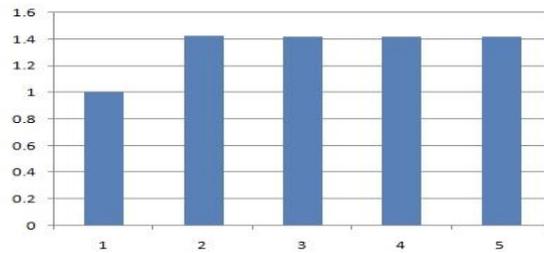

Figure 2 Convergence of $z_1 + a_1 + a_2 + a_3 + a_4$, when $x = 0.5$ and $z_1$ is assumed 1

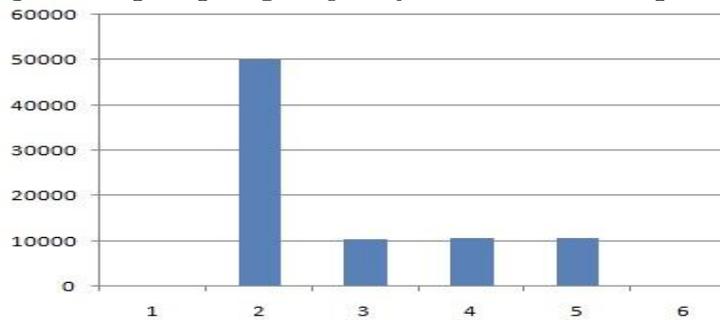

Figure 3 Convergence of $z_1 + a_1 + a_2 + a_3 + a_4$ when $x = 1$ and $z_1$ assumed is 1

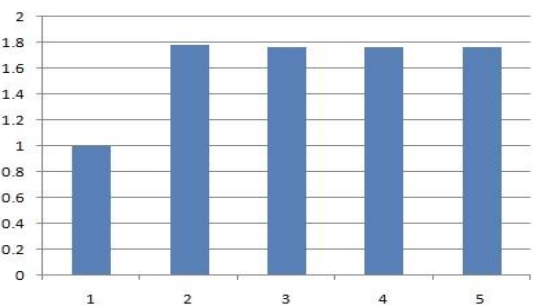

Figure 4 Convergence of $z_1 + a_1 + a_2 + a_3 + a_4$, when $x = 100$ and $z_1$ assumed is 1

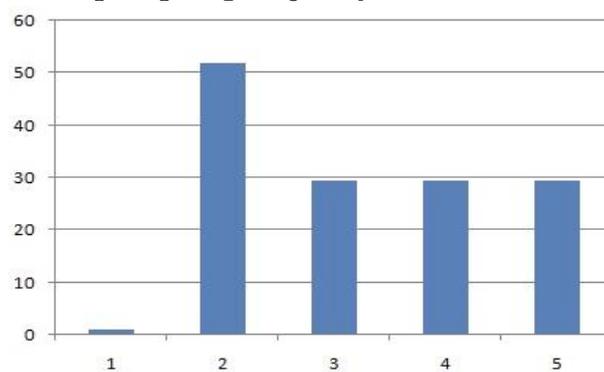

Figure 5 Convergence of $z_1 + a_1 + a_2 + a_3 + a_4$, when $x = 10^5$ and $z_1$ assumed is 1

Along y-axis, values of $z_1$, $z_2 = z_1 + a_1$, $z_3 = z_2 + a_2$, $z_4 = z_3 + a_3$ and $z_5 = z_4 + a_4$ are plotted. It is clear from these figures, the magnitudes of corrections i.e. $|z_2 - z_1|$, $|z_3 - z_2|$, $|z_4 - z_3|$ and $|z_5 - z_4|$ are in decreasing order and finally magnitude for the last correction i.e. $|z_5 - z_4|$ tends to zero.

Kindly compare heights of corrected value $z_2, z_3$ $z_4$ and $z_5$, also please note difference in height of $z_4$ with $z_5$. This difference is negligible. The figures show visual presentation of decreasing magnitude of corrections that finally give precise approximation of $z$ as $z_{n+1}$. In the Figures, number of iterations, $n = 4$, yielded precise approximation of $y$ equal to $\ln z$, where $W(x) = y$.

## 2.4. Rate Of Convergence

To find how fast the approximation reaches its actual value, rate of convergence is evaluated. Rate of convergence of the sequence $a_1, a_2, a_3, \ldots, a_n$ is found by the relation,

$$|(a_{n+1} - a_L)|/|a_n - a_L|^k,$$

where $a_L$ is the limiting value and $n$ tends to infinity. It is already explained in the foregoing sections that $a_{n+1}$ reaches value zero, when $(n+1)$ tends to infinity i. e., $a_L = 0$, therefore, $|a_{n+1} - a_L)|/|a_n - a_L|^k = |a_{n+1}|/|a_k|^k$ but before vanishing to zero at $(n+1)$, $a_n$ has negligibly small value and that results in $|a_{n+1}|/|a_n| = 0$ or rate of convergence $|(a_{n+1} - a_L)|/|a_n - a_L|^k = 0$. Thus *convergence of sequence $a_k$ is exponential* [29].

In various methods of approximation, the rate of convergence can be classified as sub-linear, linear, super-linear, quadratic, and cubic, in ascending order of efficiency. It is important to note that all these methods typically require an initial assumption close to the actual value of the variable being approximated. An initial assumption that varies within 10 to 20 percent of the actual value is acceptable. Details are given in section 8. If this variation increases, the approximation may fail to converge. This requirement ensures that the initial error (i.e., the difference between the initial assumption and the actual value) is a fraction of the actual value and significantly less than one. The following relation defines the characteristic nature of convergence,

$$E_{n+1} \propto C(E_n)^k,$$

where $E_{n+1}$ and $E_n$ stand for the error after the *(n+1)th* and *nth* iteration and $k$ is the order of convergence. This relationship holds particularly when the initial assumption is close to the actual value. However, this formula is devised for methods requiring initial guesses in the vicinity of the actual value. If the initial assumption is far from the actual value—as is the case in quadratic equation approximations—the error, defined as the difference between the initial assumption and the actual value, may be too large for the relationship to hold accurately in cases of Newton's and Halley's methods but saliently the quadratic equation approximation in spite of significant difference between initial assumption and actual value, converges *exponentially* with in three to five iterations.

The Lambert W Function involves an exponential term. Since an exponential term can not be evaluated to an exact value, the Lambert W Function can only be approximated. The precision of the result depends on how many times a correction is applied and how close the initial assumption of $z_1$ or $y_1$ is to the actual result in the first instance. If the initial assumption is near the actual value, fewer corrections are needed. From the examples provided in the tables that follow, it will be observed that after four or five corrections, precise results up to ten decimal points were obtained. If the calculation device is highly precise, applying more than five corrections can yield results accurate to beyond ten decimal points.

The value of $z_{n+1}$ has been approximated to yield precise value of $y = \ln(z_{n+1})$ where $z_{n=1}/z_n$ approximates one. Therefore, error equals to $(z - z_{n+1})/z$. Under ideal situation, $z = z_{n+1}$ when $n$ tends to infinity. Therefore, relation gives error

$$E = \frac{(a_{n+1} + a_{n+2} + a_{n+3} + \ldots + a_\infty)}{z}$$

It has already been proved that $a_1, a_2, a_3, \ldots$ are in decreasing order i.e. $|a_{n+1}| > |a_{n+2}| > |a_{n+3}| > \cdots$ so on and their magnitudes are so small that these can be neglected, therefore, percentage error for practical purposes is $100(a_{n+1}/z_{n+1})$.

## 3. Numerical Approximation of $y$, When Value of $x$ Is Given in Equation, $W_0(x) = y$ And $\infty > x \geq 0$

### 3.1. First Method

In this method, $ye^y = x$ will be written $z \ln(z) = x$, where $z = e^y$ and $z \ln(z) = x$, will be transformed to a quadratic equation. The relevant root of the equation will approximate the value of $z$, hence $y$ will be approximated as $\ln z$. Details are given hereunder.

#### 3.1.1. Quadratic Approximation of $W_0(x) = y$

Equation $W_0(x) = y$ or $y e^y = x$, on substituting $y$ with $\ln z$, results in $z \ln z = x$. Let at $z = z_1 + a$, the following equation holds good,
$$(z_1 + a_1)\ln(z_1 + a_1) = x, \tag{3.1}$$
where $z_1$ is an assumed real positive number and $a_1$ is a real number that will be found. What value should be assumed to $z_1$, will be taken up in section 3.1.2. Following the procedure as detailed in section 2, Equation (3.1) is written in quadratic form, $a_1^2 - l_1 a_1 - m_1 \simeq 0$, where
$$l_1 = -(3z_1 \ln z_1 + 2z_1 - x)/(\ln z_1 + 2), \tag{3.2}$$
$$m_1 = 2z_1(x - z_1 \ln z_1)/(\ln z_1 + 2), \tag{3.3}$$
and roots of this quadratic equation are $(1/2)\left(l_1 \pm \sqrt{l_1^2 + 4m_1}\right)$. These roots will be complex if $l^2 + 4m$ happens to be a negative quantity. Since solution of $ye^y = x$ is being determined in real domain, therefore, complex values of the roots will not be considered. If roots are found complex, then the value of $z_1$ must be changed so that values of $l_1^2 + 4m_1$, $l_2^2 + 4m_2$, $l_3^2 + 4m_3, \ldots, l_n^2 + 4m_n$ are positive quantities. It is reiterated, out of the two roots, the root,
$$a_1 = (1/2)\left(l_1 + \sqrt{l_1^2 + 4m_1}\right), \tag{3.4}$$
is considered for approximation and its other root $(1/2)\left(l_1 - \sqrt{l_1^2 + 4m_1}\right)$ is ignored. Otherwise, after $n^{th}$ correction, $z_n + a_n$, will have a negative value. Noteworthy to add, $n$ is a positive integer that corresponds to the number of times correction is applied depending upon the precision of the result needed. Thus, after first correction, $z_1 + a_1$ is rough approximation of $z$. Let $z_1 + a_1$ be equal to $z_2$, then
$$z_2 = z_1 + a_1 = z_1 + (1/2)\left(l_1 + \sqrt{l_1^2 + 4m_1}\right). \tag{3.5}$$
Approximation after second correction, will be given by the relation
$$z_3 = z_2 + a_2 = z_2 + (1/2)\left(l_2 + \sqrt{l_2^2 + 4m_2}\right), \tag{3.6}$$
where
$$l_2 = -(3z_2 \ln z_2 + 2z_2 - x)/(\ln z_2 + 2), \tag{3.7}$$
$$m_2 = 2z_2(x - z_2 \ln z_2)/(\ln z_2 + 2). \tag{3.8}$$
Proceeding in this manner, approximation after $n^{th}$ correction will be given by relation
$$z_{n+1} = z_n + a_n = z_n + (1/2)\left(l_n + \sqrt{l_n^2 + 4m_n}\right), \tag{3.9}$$
where
$$l_n = -(3z_n \ln z_n + 2z_n - x)/(\ln z_n + 2), \tag{3.10}$$
$$m_n = 2z_n(x - z_n \ln z_n)/(\ln z_n + 2). \tag{3.11}$$
In this paper, except the cases dealt in section 5, after three times correction, $(z_3 + a_3)$ yielded precise approximation of $z$ and once $(z_3 + a_3)$, was found, $y$ was determined from the relation $y = \ln(z) = \ln(z_3 + a_3)$. However, if the result demands higher accuracy, the correction can be applied $n$ times, where $n$ depends upon the extent to which the accuracy is required.

### 3.1.2. Assumption Of the Value Of $z_1$

Since $y$ can have any positive value including zero, therefore, $z = e^y$ can have minimum value 1 and maximum value infinity. Also as $ye^y = x$ or $\ln y + y = \ln x$ and if $x$ has value more than $e$, then $y = \ln x - \ln y$. Since $\ln y$ is far smaller than $y$, therefore, $y$ roughly approximates $\ln x$. In other words, $z_1$ being equal to $e^y$, roughly approximates $x$, when $x > e$. Accordingly, when $x \geq e$, assumption of value of $z_1$, equal to or less than $x$ can be made. Since $z_1$ can have minimum value one, therefore, $z_1$ *can be assumed any value between* 1 *and* $x$.

In the same way, when $x < e$, assumption of value of $z_1$ between 1 and $e$ can be made. Care will also have to be exercised, while assuming the value of $z_1$ that $z_1 + a_1$, $z_2 + a_2$, $z_3 + a_3$, ... , $z_n + a_n$ are all real and positive quantities. If root given by Equation (3.4), is found complex, then assumption of value of $z_1$ will have to be changed. Further, since $y = \ln(z_n + a_n)$, therefore, $z_n + a_n$ must be a positive real quantity. Although theory for assuming the value of $z_1$ has been detailed, *value of $z_1$ in this paper is assumed* 1 *for x, varying from* $10^{-5}$ *to* $10^5$ as is clear from the Table 3.1. But that does not mean, the value of $z_1$ should necessarily be assumed 1, it can be assumed any other value in consonance with theory already explained. Also please peruse Table 5.1, it will be seen, for given $x = 10^5$, the value of $z_1$ assumed from 1 to $10^{12}$ yielded correct results.

### 3.1.3. Veracity and Precision of Formulae Derived for approximation of y

Based upon quadratic approximation as explained in section 2, some examples are given in Table 3.1. Not only these examples prove the veracity of quadratic approximation but also reveal no error after fourth correction.

Table 3.1 Displaying accuracy of quadratic approximation for y when $W(x) = y$ is given

| S. N. | $x$ | As $z_1$ | First time corrected $z_2 = z_1 + a_1$ | Second time corrected $z_3 = z_2 + a_2$ | Third time corrected $z_4 = z_3 + a_3$ | Fourth time corrected $z_5 = z_4 + a_4$ | Calculated $y = \ln(z_5)$ | Actual $y$ | % Err |
|---|---|---|---|---|---|---|---|---|---|
| 1 | $10^{-5}$ | 1 | 1.00000999995 | Not needed | Not needed | Not needed | $9.9999(10^{-6})$ | $9.9999 \times 10^{-6}$ | 0 |
| 2 | 0.1 | 1 | 1.0956356 | 1.0955719 | Not needed | Not needed | .09127653 | .09127653 | 0 |
| 3 | 0.5 | 1 | 1.4252391 | 1.4215299 | Not needed | Not needed | .35173371 | .3517337 | 0 |
| 4 | 1.0 | 1 | 1.7807764 | 1.7632227 | 1.7632228 | Not needed | .56714329 | .56714329 | 0 |
| 5 | 100 | 1 | 51.962237 | 29.437158 | 29.536569 | 29.53659905 | 3.38563014 | 3.38563014 | 0 |
| 6 | $10^5$ | 1 | 50001.99996 | 10510.1993 | 10770.5576 | 10770.55638 | 9.2845715 | 9.2845714 | 0 |
| 7 | $10^{20}$ | 1 | $5(10^{19})$ | $2.297042(10^{18})$ | $2.36324704(10^{18})$ | $2.363688732(10^{18})$ | 42.306755092 | 42.306755092 | 0 |

It is explicit from the data given at serial numbers 1, 2 and 3 that *when $z_1 \gg a_1$, only single correction was needed and, when $z_1$ was not at all comparable with $z_1 + a_1$, three times corrections were needed.* Above detailed theory proves Lemma 3.1.

**Lemma** 3.1: *When a positive value of x between 0 and infinity is given in Lambert W function $W(x) = y$, then numerical value of y can be evaluated by quadratic approximation using relation $y = \ln(z_n + a_n)$, where $z_1$ is assumed* 1 *and n pertains to number of times the correction is applied depending upon the extent of accuracy required and*

$$z_1 + a_1 = z_1 + (1/2)\left(l_1 + \sqrt{l_1^2 + 4m_1}\right),$$
$$l_1 = -(3z_1 \ln z_1 + 2z_1 - x)/(\ln z_1 + 2),$$
$$m_1 = 2z_1(x - z_1 \ln z_1)/(\ln z_1 + 2),$$
$$z_2 + a_2 = z_2 + (1/2)\left(l_2 + \sqrt{l_2^2 + 4m_2}\right),$$
$$l_2 = -(3z_2 \ln z_2 + 2z_2 - x)/(\ln z_2 + 2),$$

$$m_2 = 2z_2(x - z_2 \ln z_2)/(\ln z_2 + 2),$$
$$...,$$
$$z_n + a_n = z_n + (1/2)\left(l_2 + \sqrt{l_n^2 + 4m_n}\right),$$
$$l_n = -(3z_n \ln z_n + 2z_n - x)/(\ln z_n + 2),$$
$$m_n = 2z_n(x - z_n \ln z_n)/(\ln z_n + 2).$$

### 3.2. Second Method

In this method, $y e^y = x$ is written as $\ln(y) + y = \ln(x)$ and, then this equation is transformed to a quadratic equation. Root of the quadratic equation will help approximating the required value of $y$. Details are given below.

### 3.1.1. Quadratic Approximation of $W_0(x) = y$

Let $y$ be $y_1 + a_1$, then
$$(y_1 + a_1) + \ln(y_1 + a_1) = \ln(x). \tag{3.12}$$
As explained in section 2, this equation is transformed to a quadratic equation,
$$a_1^2 - l_1 a - m_1 = 0, \tag{3.13}$$
where
$$l_1 = -\{3y_1 + 2 - \ln(x/y_1)\}, \tag{3.14}$$
$$m_1 = -2y_1\{y_1 - \ln(x/y_1)\}. \tag{3.15}$$
Out of the two roots of the quadratic Equation (3.13), root
$$a_1 = (1/2)\left(l_1 + \sqrt{l_1^2 + 4m_1}\right), \tag{3.16}$$
is taken into consideration and its other root $(1/2)\left(l_1 - \sqrt{l_1^2 + 4m_1}\right)$ is neglected on account of the fact that it leads to a value $z_n + a_n$ which is a negative quantity. Thus, after first correction, $y_1 + a_1$ is the rough approximation of $y$ and if assumption of $y_1$ happens to be such that $y_1 \gg a_1$, then even after first correction, result will be accurate. Kindly refer to Table 3.2.

Let $y_1 + a_1$ be equal to $y_2$, then
$$y_2 = y_1 + a_1 = y_1 + (1/2)\left(l_1 + \sqrt{l_1^2 + 4m_1}\right). \tag{3.17}$$
Approximation after second correction will be given by relation
$$y_3 = y_2 + a_2 = y_2 + (1/2)\left(l_2 + \sqrt{l_2^2 + 4m_2}\right), \tag{3.18}$$
where
$$l_2 = -\{3y_2 + 2 - \ln(x/y_2)\}, \tag{3.19}$$
$$m_2 = -2y_2\{y_2 - \ln(x/y_2)\}. \tag{3.20}$$
Proceeding in this manner, approximation after n$^{th}$ correction will be given by the relation
$$y_{n+1} = y_n + a_n = y_n + (1/2)(l_n + \sqrt{l_n^2 + 4m_n}), \tag{3.21}$$
where
$$l_n = -\{3y_n + 2 - \ln(x/y_n)\}, \tag{3.22}$$
$$m_n = -2y_n\{y_2 - \ln(x/y_n)\}. \tag{3.23}$$
In this paper, correction is applied three times and $(y_3 + a_3)$ yielded the precise approximation. For determining rate of convergence, explanation given in section 2.4 is applicable and is not repeated for the sake of brevity.

### 3.2.2. Assumption Of the Value Of $y_1$

Since $y_1$ can have any positive value and $y = \ln x - \ln y$, therefore, it is true that $1 \leq y < \ln(x)$, when $x \geq e$ and $y_1$ can be assigned any value between 1 and $\ln x$. When $0 \leq x < e$, then it is true that $\ln(x) \geq \ln(y)$ or $0 < y \leq x$, accordingly assumption of value of $y_1$ is made between 0 and $x$. Care will also have to be exercised while assuming the value of $y_1$, so that $y_1 + a_1, y_2 + a_2, y_3 + a_3, ..., y_n + a_n$ are all real and positive quantities, if these happen to be negative or complex, then assumption of $y_1$ will have to be changed. Although theory for assuming the value of $y_1$ has been detailed above, value of $y_1$ in this paper is assumed as 1 when $x \geq 1$ and equal to $x$ when $0 < x < 1$. Kindly peruse Table 3.2.

Table 3.2 Displaying accuracy of quadratic approximation for y when $W(x) = y$ is given

| S. N. | Given $x$ | Assumed $y_1$ | First time corrected $y_1 + a_1$ | Second time corrected $y_1 + a_1 + a_2$ | Third time corrected $y = y_1 + a_1 + a_2 + a_3$ | Actual $y$ | % age Error After Correction |
|---|---|---|---|---|---|---|---|
| 1 | $10^{-5}$ | $10^{-5}$ | $9.9999(10^{-6})$ | Not needed | Not needed | $9.9999(10^{-6})$ | 0.000000 |
| 2 | $10^{-2}$ | $10^{-2}$ | .00990147305 | .00990147384 | Not needed | .00990147384 | 0.00000 |
| 3 | $10^{-1}$ | $10^{-1}$ | .09127122105 | .09127653616 | .0912765271 | .09127652716 | 0.00000 |
| 4 | 1 | 1 | .5615528128 | .5671433197 | .56714329 | .56714329 | 0.000000 |
| 5 | $10^2$ | 1 | 3.49503992 | 3.385628701 | 3.38563014 | 3.3856301403 | $8.86098(10^{-9})$ |
| 6 | $10^5$ | 1 | 9.880553811 | 9.28455331 | 9.284571428 | 9.2845714286 | $6.46233(10^{-9})$ |
| 7 | $10^{10}$ | 1 | 21.20598257 | 20.02867132 | 20.02868541 | 20.028685413 | $1.49785(10^{-8})$ |
| 8 | $10^{20}$ | 1 | 44.14031445 | 42.3067489 | 42.30675509 | 42.306755096 | $1.418213(10^{-8})$ |
| 9 | $10^{50}$ | 1 | 113.16429118 | 110.4249176 | 110.4249188 | 110.42491883 | $2.71678(10^{-8})$ |
| 10 | $10^{500}$ | 1 | 688.7813268 | 684.2625008 | 684.2472086 | 684.2472086 | $4.384380(10^{-9})$ |

It is explicit from the data given at serial numbers 1 and 2 that correction was needed once or twice for approximating value of $y$. Above detailed theory proves Lemma 3.2.

**Lemma** 3.2: *When a positive value of $x$ is given in Lambert W function $W(x) = y$, then numerical value of $y$ can be evaluated by the quadratic approximation using relation $y = y_n + a_n$ where $y_1$ is assumed as 1, when $x \geq 1$; as equal to $x$, when $0 < x < 1$; n pertains to number of times the correction is applied depending upon the extent of accuracy required and*

$$y_1 + a_1 = y_1 + (1/2)\left(l_1 + \sqrt{l_1^2 + 4m_1}\right),$$
$$l_1 = -\{3y_1 + 2 - \ln(x/y_1)\},$$
$$m_1 = -2y_1\{y_1 - \ln(x/y_1)\},$$
$$y_2 + a_2 = y_2 + (1/2)\left(l_2 + \sqrt{l_2^2 + 4m_2}\right),$$
$$l_2 = -\{3y_2 + 2 - \ln(x/y_2)\},$$
$$m_2 = -2y_2\{y_2 - \ln(x/y_2)\}.$$
$$...,$$
$$y_n + a_n = y_n + (1/2)\left(l_n + \sqrt{l_n^2 + 4m_n}\right),$$
$$l_n = -\{3y_n + 2 - \ln(x/y_n)\},$$
$$m_n = -2y_n\{y_2 - \ln(x/y_n)\}.$$

## 4. Numerical Approximation Of $y$, When Value Of $-x$ Is Given In Equation $W(-x) = -y$, Where $0 \leq x \leq 1/e$

### 4.1. First Method

In this method, $-ye^{-y} = -x$ is written as $(1/z)\ln(z) = x$, where $z = e^y$ and $(1/z)\ln(z) = x$ is transformed to to a quadratic equation. The relevant root of the equation will help approximating the value of $z$ and that will ease approximation of $y$.

#### 4.1.1. Quadratic Approximation of $W(-x) = -y$

Equation $W(-x) = -y$ or $-y e^{-y} = -x$, on substituting $y$ with $\ln z$, results in $(1/z)\ln z = x$, where $x$ is a positive real quantity excluding infinity. Let $z$ be equal to $z_1 + a_1$, where $z$, and $z_1$ are real positive quantities, then

$$\ln(z_1 + a_1) = x(z_1 + a_1). \tag{4.1}$$

Term $ln(z_1 + a_1)$ can be written $\ln z_1 + \ln(1 + a_1/z_1)$. Using Equation (2.3), term $ln(z_1 + a_1)$ can further be written $\ln z_1 + 2a_1/(2z_1 + a_1)$. Putting this value for $ln(z_1 + a_1)$ in Equation (4.1), it takes the form,
$$\ln z_1 + 2a_1/(2z_1 + a_1) - x(z_1 + a_1) = 0.$$
On simplification, it, yields a quadratic equation $a_1^2 - l_1 a - m_1 = 0$ in variable $a_1$, which has roots $(1/2)(l_1 \pm \sqrt{l_1^2 + 4m_1})$, where
$$l_1 = -(3z_1 x - \ln z_1 - 2)/x,$$
$$m_1 = 2z_1 (\ln z_1 - z_1 x)/x.$$

Therefore, after first correction,
$$z_2 = z_1 + a_1 = z_1 + (1/2)(l_1 + \sqrt{l_1^2 + 4m_1}),$$
after second correction,
$$z_3 = z_2 + a_2 = z_2 + (1/2)(l_2 + \sqrt{l_2^2 + 4m_2}),$$
$$l_2 = -(3z_2 x - \ln z_2 - 2)/x,$$
$$m_2 = 2z_2 (\ln z_2 - z_2 x)/x.$$
...,

after $n^{th}$ correction,
$$z_{n+1} = z_n + a_n = z_n + (1/2)(l_n + \sqrt{l_n^2 + 4m_n}),$$
$$l_n = -(3z_n x - \ln z_n - 2)/x,$$
$$m_n = 2z_n (\ln z_n - z_n x)/x.$$

Earlier while devising methods for numerical approximation of $y$ in Lambert W function $W_0(x) = y$, where value of $x$ was given, one root, out of two roots of the quadratic equation was considered on account of the fact that the other root yielded negative value of $z_n + a_n$ leading to a complex value of $y$ and hence was rejected. But in the present case, *both roots $a_1 = (1/2)(l_1 + \sqrt{l_1^2 + 4m_1})$ and $a_1' = (1/2)(l_1 \pm \sqrt{l_1^2 + 4m_1})$, give positive values of $z_n + a_n$, $z_n' + a_n'$, hence both are being considered* as is explicit from Table 4.1.

For second value of $z$, other root will be denoted, $a_1' = (1/2)(l_1 - \sqrt{l_1^2 + 4m_1})$. Following the same procedure, after first correction,
$$z_2' = z_1 + a_1' = z_1 + (1/2)(l_1 - \sqrt{l_1^2 + 4m_1}),$$
after second correction,
$$z_3' = z_2' + a_2' = z_2' + (1/2)(l_2' - \sqrt{l_2'^2 + 4m_2'}),$$
$$l_2' = -(3z_2' x - \ln z_2' - 2)/x,$$
$$m_2' = 2z_2' (\ln z_2' - z_2' x)/x.$$
...,

after $n^{th}$ correction,
$$z_{n+1}' = z_n' + a_n' = z_n' + (1/2)(l_n' - \sqrt{l_n'^2 + 4m_n'}),$$
$$l_n' = -(3z_n' x - \ln z_n' - 2)/x,$$
$$m_n' = 2z_n' (\ln z_n' - z_n' x)/x.$$

### 4.1.2. Assumption Of the Value Of $z_1$

In equation, $-y e^{-y} = -x$, when $y$ varies from 0 to minus infinity, $x$ takes the value from 0 to 0 but it is not a flat straight line, $x = 0$. When $y$ rises above 0, value of $-y e^{-y}$ decreases and reaches minimum value and then increases and reaches 0. To find out the minimum value of $-x$, derivative of $-ye^{-y}$ is equated with zero and this yields $-y = -1$ or $-x = -1/e$. When $y = 0$, then $x = 0$, when $-y = -1$, then $-x = -1/e$. In such situation $z = e^y$ will vary from 1 to $e$. When $-y = -1$, then $-x = -1/e$, when $-y = -\infty$, then $x = 0$. In such situation $z = e^y$ will vary from $e$ to $\infty$. Therefore, one value of $x$ will correspond to values of $y$. Again value of $z$ should be such that $z_n + a_n$, $z_n' + a_n'$ must be positive real number. Referring to Table 5.2 in section 5, in equation $W(-x) = -y$, for value of $-x$ assumed $-.1$, $z_1$ was assigned values 0.3, $10^3$ and $10^9$ and such assignments yielded precise approximation of $y = \ln(z_4 + a_4)$. If the calculator carrying out the computations, happens to be precise, $z_1$ can be assigned value more than $10^9$. In equation $W(-x) = -y$, for value of $-x$ assumed $-.1$, $z_1$ was assigned values .3 and 5 and such assignments yielded precise approximation of second value of $y = \ln(z_4' + a_4')$. However, in Table 4.1, value of $z_1$ is assumed 2 and that facilitated approximation of both values of $y$.

Table 4.1 Displaying accuracy of quadratic approximation for y when x is given in Equation $W(-x) = -y$

| S. N. | $x$ | As s. $z_1$ | First time corrected $z_1 + a_1$ | Second time corrected $z_1 + a_1 + a_2$ | Third time corrected $z_1 + a_1 + a_2 + a_3$ | Calculated $y = \ln(z_1 + a_1 + a_2 + a_3)$ | Actual $y$ | Percentage error after correction |
|---|---|---|---|---|---|---|---|---|
| 1a | .365 | 2 | 2.956051512 | 3.097484097 | 3.09768805 | 1.130656043 | 1.1306553125 | $6.4608(10^{-5})$ |
| 1b | .365 | 2 | 2.422430039 | 2.410565801 | 2.410465598 | .879819923 | .879819986 | $7.1605(10^{-6})$ |
| 2a | .25 | 2 | 7.35020321 | 8.610707527 | 8.613169456 | 2.153292364 | 2.153292364 | 0.00000 |
| 2b | .25 | 2 | 1.422385512 | 1.429611849 | 1.429611805 | .3574029424 | .3574029562 | $3.8612(10^{-6})$ |
| 3a | .10 | 2 | 23.83488834 | 35.94782488 | 35.77152074 | 3.577152067 | 110.42491883 | $2.7168(10^{-8})$ |
| 3b | .10 | 2 | 1.096588361 | 1.118326385 | 1.118325592 | .1118325595 | .1118325592 | $-2.5354(10^{-7})$ |
| 4a | $10^{-3}$ | 2 | 2689.157469 | 8974.76983 | 9118.006099 | 9.11800643 | 9.11800647 | $4.3870(10^{-7})$ |
| 4b | $10^{-3}$ | 2 | .971574359 | 1.001003721 | 1.001001503 | .0010010018 | .0010010015 | $2.6667(10^{-5})$ |

It is pertinent to mention that serial numbers 1a, 2a, 3a and 4a correspond to data relating to root $(1/2)(l_1 + \sqrt{l_1^2 + 4m_1})$ and serial numbers 1b, 2b, 3b and 4b correspond to data relating to root $(1/2)(l_1 - \sqrt{l_1^2 + 4m_1})$ for both cases. *It is not necessary that value of $z_1$ should be assumed same for both roots, value of $z_1$ pertaining to root $(1/2)(l_1 + \sqrt{l_1^2 + 4m_1})$ can be assumed different pertaining to root $(1/2)(l_1 - \sqrt{l_1^2 + 4m_1})$. That proves Lemma 4.1.*

**Lemma** 4.1: *When given value of $x$ is such that $x \leq 1/e$ in $-ye^{-y} = -x$ or Lambert W function $W(-x) = -y$, then one numerical value of $y$ can be evaluated by the quadratic approximation and using relation $y = \ln(z_n + a_n)$, where*

$$z_2 = z_1 + a_1 = z_1 + (1/2)(l_1 + \sqrt{l_1^2 + 4m_1}),$$
$$l_1 = -(3z_1 x - \ln z_1 - 2)/x,$$
$$m_1 = 2z_1 (\ln z_1 - z_1 x)/x,$$
$$z_3 = z_2 + a_2 = z_2 + (1/2)(l_2 + \sqrt{l_2^2 + 4m_2}),$$
$$l_2 = -(3z_2 x - \ln z_2 - 2)/x,$$
$$m_2 = 2z_2 (\ln z_2 - z_2 x)/x.$$
$$\ldots,$$
$$z_{n+1} = z_n + a_n = z_n + (1/2)(l_n + \sqrt{l_n^2 + 4m_n}),$$
$$l_n = -(3z_n x - \ln z_n - 2)/x,$$
$$m_n = 2z_n (\ln z_n - z_n x)/x.$$

*and second numerical value of $y$ can be evaluated by quadratic approximation and using relation $y = \ln(z_n' + a_n')$ where*

$$z_2' = z_1 + a_1' = z_1 + (1/2)(l_1 - \sqrt{l_1^2 + 4m_1}),$$
$$z_3' = z_2' + a_2' = z_2' + (1/2)(l_2' - \sqrt{l_2'^2 + 4m_2'}),$$
$$l_2' = -(3z_2' x - \ln z_2' - 2)/x,$$
$$m_2' = 2z_2' (\ln z_2' - z_2' x)/x.$$

$$z'_{n+1} = z'_n + a'_n = z'_n + (1/2)(l'_n - \sqrt{l'^2_n + 4m'_n}),$$
$$l'_n = -(3z'_n x - \ln z'_n - 2)/x,$$
$$m'_n = 2z'_n (\ln z'_n - z'_n x)/x.$$

$z_1$ is assumed as 2, and $n$ pertains to number of times the correction is applied depending upon the extent of accuracy required.

### 4.2. Second Method

In this method, $-ye^{-y} = -x$ will be written as $\ln(y_1 + a_1) - y_1 - a_1 = \ln(x)$ and, then this equation is transformed to a quadratic equation in $a_1$. Root of quadratic equation will provide correction to $y_1$.

#### 4.2.1. Quadratic Approximation of $W(-x) = -y$

Let at $y = y_1 + a_1$, the following equation holds good,
$$-(y_1 + a_1) + \ln(y_1 + a_1) = \ln(x).$$
This equation can be written
$$-a + \ln(1 + a_1/y_1) + \ln(y_1) = \ln(x) + y_1.$$
Using Equation (2.3) and simplification yields
$$-a_1 + 2/(2y_1/a_1 + 1) - \ln(x/y_1) - y_1 = 0.$$
which is quadratic equation $a_1^2 - l_1 a_1 - m_1 = 0$, where
$$l_1 = -\{3y_1 - 2 + \ln(x/y_1)\},$$
$$m_1 = -2y_1\{y_1 + \ln(x/y_1)\}.$$
Roots of this equation are $(1/2)(l \pm \sqrt{l^2 + 4m})$ and both roots will give corrections as explained in section 4.1.1. Therefore,
$$y_2 = y_1 + a_1 = y_1 + (1/2)\left(l_1 + \sqrt{l_1^2 + 4m_1}\right),$$
$$y_3 = y_2 + a_2 = y_2 + (1/2)\left(l_2 + \sqrt{l_2^2 + 4m_2}\right),$$
$$l_2 = -\{3y_2 - 2 + \ln(x/y_2)\},$$
$$m_2 = -2y_2\{y_2 + \ln(x/y_2)\},$$
$$...,$$
$$y_{n+1} = y_n + a_n = y_n + (1/2)(l_n + \sqrt{l_n^2 + 4m_n}),$$
$$l_n = -\{3y_n - 2 + \ln(x/y_n)\},$$
$$m_n = -2y_n\{y_n + \ln(x/y_n)\}.$$
Similarly for second value $y'_n + a'_n$,
$$y'_2 = y_1 + a'_1 = y_1 + (1/2)\left(l_1 - \sqrt{l_1^2 + 4m_1}\right),$$
$$y'_3 = y'_2 + a'_2 = y'_2 + (1/2)\left(l'_2 - \sqrt{l'^2_2 + 4m'_2}\right),$$
$$l'_2 = -\{3y'_2 - 2 + \ln(x/y'_2)\},$$
$$m'_2 = -2y'_2\{y'_2 + \ln(x/y'_2)\},$$
$$...,$$
$$y'_{n+1} = y'_n + a'_n = y'_n + (1/2)\left(l'_n - \sqrt{l'^2_n + 4m'_n}\right),$$
$$l'_n = -\{3y'_n - 2 + \ln(x/y'_n)\},$$
$$m'_n = -2y'_n\{y'_n + \ln(x/y'_n)\}.$$

#### 4.2.2. Assumption Of the Value Of $y_1$

From table 4.2, it is clear, for a given value of $-x$, Lambert W function $W(-x) = -y$ has two corresponding values of $-y$, one value $-y$ lies between $-1$ and $0$ whereas other between $-1$ and $-\infty$ and smaller value of $-y$ is such that $y \geq x$ or $-y \leq -x$ and larger value of $y$ is such that $y > 0$ or $-y < 0$. That means for smaller value of $y$, assumption of $y_1$ can be made as if it were equal to or more than $x$ and for larger value of $y$, assumption can be made as if $y_1$ were more than zero.

Care will also have to be exercised, while assuming the value of $y_1$ that roots $y_1 + a_1$, $y_2 + a_2$, $y_3 + a_3, \ldots, y_n + a_n$ and $y_1 + a_1'$, $y_2' + a_2'$, $y_3' + a_3', \ldots, y_n' + a_n'$ are all real and positive. If a root is found complex, then assumption of $y_1$ will have to be changed. Since numerical approximation of $y$ in real domain is being found, that requires assumption of $y_1$ must be a real quantity such that it should result in value of roots $a_1, a_2, a_3 \ldots, a_n$ and $a_1', a_2', a_3' \ldots, a_n'$ in real domain. To prove veracity of formulae derived for approximation of $y$, some examples are given in Table 4.2.

Table 4.2 Displaying accuracy of quadratic approximation for y when x is given in Equation $W(-x) = -y$

| S. N. | Given $x$ | Ass. $y_1$ | First time corrected $y_1 + a_1$ | 2nd time corrected $y_1 + a_1 + a_2$ | 3rd time corrected $y = y_1 + a_1 + a_2 + a_3$ | Actual $y$ | % age error after correction |
|---|---|---|---|---|---|---|---|
| 1a | .365 | 1 | 1.129352923 | 1.130655313 | 1.130655313 | 1.130655313 | 0.0000 |
| 1b | .365 | 1 | 0.7570090661 | 0.882015716 | 0.879820082 | 0.879820092 | $1.1365(10^{-6})$ |
| 2a | 0.1 | 1 | 3.39179597 | 3.57713465 | 3.577152064 | 3.577152064 | 0.0000 |
| 2b | 0.1 | 0.1 | 0.1118472693 | 0.1118325592 | No need | 0.1118325592 | 0.0000 |
| 3a | $10^{-3}$ | $10^{-3}$ | 8.486085012 | 9.117971815 | 9.118006532 | 9.11800647 | $6.7997(10^{-7})$ |
| 3b | $10^{-3}$ | $10^{-3}$ | 0.00100100150 | No need | No need | 0.00100100150 | 0.00000 |

In the Table 4.2, abbreviation 'Ass.' denotes 'Assumed.' In some columns, it is mentioned, 'No need,' that says, there is no necessity of applying correction. At serial 1b and 2a, small error in $10^{-6}$ is noted. This error can also be eliminated by applying correction fourth time. Precision of formulae derived above, proves Lemma 4.2a and 4.2b.

**Lemma** 4.2a: *When given value of $x$ is such that $x \leq 1/e$ in $-ye^{-y} = -x$ or Lambert W function $W(-x) = -y$, then one of numerical value of $y$ can be evaluated by quadratic approximation and using relation $y = y_n + a_n$ where $y_1$ is assumed $1$, $n$ pertains to number of times the correction is applied depending upon the extent of accuracy required,*

$$y_2 = y_1 + a_1 = y_1 + (1/2)\left(l_1 + \sqrt{l_1^2 + 4m_1}\right),$$
$$l_1 = -\{3y_1 - 2 + \ln(x/y_1)\},$$
$$m_1 = -2y_1\{y_1 + \ln(x/y_1)\}.$$
$$y_3 = y_2 + a_2 = y_2 + (1/2)\left(l_2 + \sqrt{l_2^2 + 4m_2}\right),$$
$$l_2 = -\{3y_2 - 2 + \ln(x/y_2)\},$$
$$m_2 = -2y_2\{y_2 + \ln(x/y_2)\},$$
$$\ldots,$$
$$y_{n+1} = y_n + a_n = y_n + (1/2)\left(l_n + \sqrt{l_n^2 + 4m_n}\right),$$
$$l_n = -\{3y_n - 2 + \ln(x/y_n)\},$$
$$m_n = -2y_n\{y_n + \ln(x/y_n)\}.$$

**Lemma** 4.2b: *When a negative value of $-x > -1/e$ is given in Lambert W function $W(-x) = -y$, then other numerical value of $y$ can be evaluated by quadratic approximation and using relation $y = y_n' + a_n'$ where $y_1$ is assumed equal to $x$, $n$ pertains to number of times the correction is applied depending upon the extent of accuracy required,*

$$y_2' = y_1 + a_1' = y_1 + (1/2)\left(l_1 - \sqrt{l_1^2 + 4m_1}\right),$$
$$l_1 = -\{3y_1 - 2 + \ln(x/y_1)\},$$
$$m_1 = -2y_1\{y_1 + \ln(x/y_1)\}.$$
$$y_3' = y_2' + a_2' = y_2' + (1/2)\left(l_2' - \sqrt{l_2'^2 + 4m_2'}\right),$$

$$l_2' = -\{3y_2' - 2 + ln(x/y_2')\},$$
$$m_2' = -2y_2'\{y_2' + ln(x/y_2')\},$$
$$...,$$
$$y_{n+1}' = y_n' + a_n' = y_n' + (1/2)\left(l_n' - \sqrt{l_n'^2 + 4m_n'}\right),$$
$$l_n' = -\{3y_n' - 2 + ln(x/y_n')\},$$
$$m_n' = -2y_n'\{y_n' + ln(x/y_n')\}.$$

## 5. Wide Range of Initial Assumptions

Although this method gives details how to make initial assumptions in sections 3 and 4, it is practically observed that in addition, a broad range of values can be considered for initial assumptions. Following sections 5.1) and 5.2 deal with such assumptions.

### 5.1 Assignment of Different Values of $z_1$ Yields Same Solution Of $W(x) = y$

In this section, it will be explained, how in equation $W(x) = y$, assumptions of different values of $z_1$ leads to the same precise result of $y$, where $y = ln\, z$ and $z = z_n + a_n$. In the Table 5.1, equation $W(10^5) = y$, is solved by assuming $z_1$ equal to $1, 10, 10^2, 10^3, 10^4, 10^5, 10^6, 10^{12}$ where $z_2, z_3, ... z_{n+1}, l_1, l_2,, ..., l_n, m_1, m_2, ..., m_n$ are given by Equations (3.5), (3.6),..., (3.9), (3.2), (3.7),..., (3.10), (3.3), (3.8),..., (3.10) respectively. It was also observed, when $x > 10^{12}$, value of $z_2$ which is $z_1 + a_1$ calculates zero, using Desmos Scientific calculator and natural logarithm of 0 is not defined, hence subsequent values of $z_3, z_4, ..., z_{n+1}$ can not be calculated, therefore, assumption of value of $z_1$ can not be made more than $10^{12}$. If the calculator happens to be precise, one can go beyond $10^{12}$. It is also noteworthy that $z_1$ can be assigned any value between 1 and $10^{12}$, when using *Desmos* Scientific calculator and the value of $y = ln(z_{n+1})$ reaches precise value though number of steps varies. This property adds playfulness and recreational elements to this method.

Table 5.1 Displaying assumptions of different values of $z_1$ leading to same results

| $x$ | Ass. $z_1$ | First time corrected $z_2$ | Second time corrected $z_3$ | Third time corrected $z_4$ | Fourth time corrected $z_5$ | Calculated $y = ln\,(z_5)$ | Actual $y$ | Err. |
|---|---|---|---|---|---|---|---|---|
| $10^5$ | 1 | 50001.99996 | 10510.1993 | 10770.5576 | 10770.55638 | 9.284571429 | 9.284571429 | 0.00 |
| $10^5$ | 10 | 49574.91369 | 10573.99124 | 10770.55692 | 10770.55638 | 9.284571429 | 9.284571429 | 0.00 |
| $10^5$ | $10^2$ | 15199.81914 | 10767.02713 | 10770.55638 | No need | 9.284571429 | 9.284571429 | 0.00 |
| $10^5$ | $10^3$ | 11639.69298 | 10770.51561 | 10770.55638 | No need | 9.284571429 | 9.284571429 | 0.00 |
| $10^5$ | $10^4$ | 10770.59206 | 10770.55638 | No need | No need | 9.284571429 | 9.284571429 | 0.00 |
| $10^5$ | $10^5$ | 10120.89002 | 10770.57739 | 10770.55638 | No need | 9.284571429 | 9.284571429 | 0.00 |
| $10^5$ | $10^6$ | 8439.5434 | 10771.81686 | 10770.55638 | No need | 9.284571429 | 9.284571429 | 0.00 |
| $10^5$ | $10^{12}$ | 12000 | 10770.44626 | 10770.55638 | No need | 9.284571429 | 9.284571429 | 0.00 |

### 5.2. Assignment Of Different Values of $z_1$ Yields Same Solution Of $W(-x) = -y$

Now it will be shown, how in equation $W(-x) = -y$ or $-ye^{-y} = -x$, assumptions of different values of $z_1$ leads to the same precise results of $y$, where $y = ln\, z$ and $z = z_{n+1} = z_n + a_n$. It has already been proved that solution of equation $-ye^{-y} = -x$ yields two values of $y$ for a given single value of $x$, such that $0 < x \leq 1/e$ and one value of $-y$ lies between 0 and $-1$ and other between $-1$ and minus infinity. In Table 5.2, equation $W(-0.1) = -y$, is solved by assuming $z_1$ equal to $0.3, 10^3, 10^9$, where values of $z_2, z_3, ..., z_{n+1}, l_1, l_2,, ..., l_n, m_1, m_2, ..., m_n$ are given in section 4.1. Referring to Table 5.2, it was also observed that when $x > 10^9$, value of $z_2$ which is $z_1 + a_1$ calculates zero, using *Desmos* Scientific calculator and natural logarithm of 0 is not defined, hence subsequent values of $z_3, z_4, ..., z_{n+1}$ can not be calculated, therefore, assumption of value of $z_1$ can not be made more than $10^9$. If the calculator happens to be precise, one can go beyond $10^9$.

Table 5.2 Displaying assumptions of different values of $z_1$ leading to same results for $-ye^{-y} = -x$ and $-y > -1$

| $x$ | Ass. $z_1$ | First time corrected $z_2$ | Second time corrected $z_3$ | Third time corrected $z_4$ | Fourth time corrected $z_5$ | Calculated $y = ln\,(z_5)$ | Actual $y$ | Error |
|---|---|---|---|---|---|---|---|---|
| 0.1 | 0.3 | 6.079140914 | 31.59402227 | 35.76931211 | 35.77152064 | 3.577152064 | 3.577152064 | Nil |
| 0.1 | $10^3$ | 51.019259 | 35.82203151 | 35.77152064 | No need | 3.577152064 | 3.577152064 | Nil |
| 0.1 | $10^9$ | 180 | 39.0621895 | 35.77230776 | 35.77152064 | 3.577152064 | 3.577152064 | Nil |

Table 5.3 Displaying assumptions of different values of $z_1$ leading to same results for $-ye^{-y} = -x$ and $0 < -y < -1$

| $x$ | Ass. $z_1$ | First time corrected $z_2$ | Second time corrected $z_3$ | Third time corrected $z_4$ | Fourth time corrected $z_5$ | Calculated $y = \ln(z_5)$ | Actual $y$ |
|---|---|---|---|---|---|---|---|
| 0.1 | 0.3 | 1.58113103 | 1.113880206 | 1.118325962 | 1.118325598 | 0.1118325649 | 0.1118325592 |
| 0.1 | 5 | 0.641251085 | 1.137654313 | 1.118325064 | 1.118325591 | 0.1118325586 | 0.1118325592 |

Using formulae for $z_2, z_3, \ldots z_{n+1}$, $l_1, l_2, \ldots, l_n$, $m_1, m_2, \ldots, m_n$ given in section 4.1, the equation $W(-.1) = -y$, is solved by assuming $z_1$ equal to 0.3 and 5.0. These solutions are entered in Table 5.3. The solutions depicted in Tables 5.1 and 5.2 prove different assignment of different values of $z_1$ yield same results.

## 6. Algorithm For Solving Lambert W Function $W(x) = Y$ And Demonstration Of Its Use for Solving $W(10^{20}) = y$

1. Check whether real value of $x \geq 0$ or $x < 0$. If $x < 0$, go to 15.
2. If $x = 0$, go to 35.
3. Let $z_1 = 1$.
4. Compute
$$l_1 = -(3z_1 \ln z_1 + 2z_1 - x)/(\ln z_1 + 2)$$
$$= -\{3(1)\ln(1) + 2(1) - 10^{20}\}/\{\ln(1) + 2\} = 5(10^{19})$$
$$m_1 = 2z_1(x - z_1 \ln z_1)/(\ln z_1 + 2)$$
$$= 2(1)\{10^{20} - (1)\ln(1)\}/\{\ln(1) + 2\} = 5(10^{19})\} = 1(10^{20})$$
5. Compute
$$z_2 = z_1 + (1/2)(l_1 + \sqrt{l_1^2 + 4m_1})$$
$$= 1 + (1/2)[5(10^{19}) + \{(5^2)(10^{19})^2 + 4(1)(10^{20})\}^{1/2} = 5(10^{19})$$
6. $l_2 = -(3z_2 \ln z_2 + 2z_2 - x)/(\ln z_2 + 2)$,
$$= -[3(5)(10^{19})\ln\{5(10^{19})\} + 2(5)(10^{19}) - 10^{20}]/[\ln\{(5)(10^{19})\} + 2]$$
$$= -1.43665347(10^{20})$$
$$m_2 = 2z_2(x - z_2 \ln z_2)/(\ln z_2 + 2)$$
$$= 2(5)(10^{19})[10^{20} - (5)(10^{19})\ln\{(5)(10^{19})\}]/[\ln\{(5)(10^{19})\} + 2]$$
$$= -4.57768981(10^{39})$$
7. Compute
$$z_3 = z_2 + (1/2)(l_2 + \sqrt{l_2^2 + 4m_2})$$
$$= 5(10^{19}) + (1/2)[-(1.43665347(10^{20})$$
$$+ \{(-1.43665347)^2(10^{19})^2 + 4(-4.57768982)(10^{39})\}^{1/2}$$
$$= 2.297042(10^{18})$$
8. Compute
$$l_3 = -(3z_3 \ln z_3 + 2z_3 - x)/(\ln z_3 + 2)$$
$$= -[3(2.297042)(10^{18})\ln\{(2.297042)(10^{18})\} + 2(2.297042)(10^{18})$$
$$- 10^{20}]/[\ln\{(2.297042)(10^{18}\} + 2] = -4.45558155(10^{18})$$
$$m_3 = 2z_3(x - z_3 \ln z_3)/(\ln z_3 + 2)$$
$$= 2(2.297042)(10^{18})[10^{20} - (2.297042)(10^{18})\ln\{(2.297042)(10^{18})\}]$$
$$/[\ln\{(2.297042)(10^{18})\} + 2] = 2.99365075(10^{35})$$
9. Compute
$$z_4 = z_3 + (1/2)(l_3 + \sqrt{l_3^2 + 4m_3})$$
$$= 2.297042(10^{18}) + (1/2)\{(-4.45558155)^2(10^{18})^2 + 4(2.99365075)(10^{35})\}^{1/2}$$
$$= 2.36324704(10^{18})$$
10. Compute
$$l_4 = -(3z_4 \ln z_4 + 2z_4 - x)/(\ln z_4 + 2)$$
$$= -[3(2.36324704)(10^{18})\ln\{(2.36324704)(10^{18})\} + 2(2.36324704)(10^{18})$$
$$- 10^{20}]/[\ln\{(2.36324704)(10^{18}\} + 2] = -4.61938531(10^{18})$$

$$m_4 = 2z_4(x - z_4 \ln z_4)/(\ln z_4 + 2)$$
$$= 2(2.36324704)(10^{18})[10^{20} - (2.36324704)(10^{18})\ln\{(2.36324704)(10^{18})\}]$$
$$/[\ln\{(2.36324704)(10^{18})\} + 2] = \mathbf{2.0405081(10^{33})}$$

11. Compute
$$z_5 = z_4 + (1/2)\left(l_4 + \sqrt{l_4^2 + 4m_4}\right)$$
$$= 2.36324704(10^{18})$$
$$+ (1/2)\{(-4.61938531)^2(10^{18})^2 + 4(2.0405081)(10^{35})\}^{1/2}$$
$$= \mathbf{2.363688732(10^{18})}$$

12. Compute $\ln(z_5)$. Let it be $Y = \mathbf{42.306755092}$
13. Go to 34.
14. Let $X = -x$, check if $\infty > X > 1/e$. If yes go to 37.
15. Let $z_1 = 2$.
16. Compute
$$l_1 = -(3z_1 X - \ln z_1 - 2)/X,$$
$$m_1 = 2z_1 (\ln z_1 - z_1 X)/X.$$
18. Compute
$$z_2 = z_1 = z_1 + (1/2)(l_1 + \sqrt{l_1^2 + 4m_1}).$$
19. Compute
$$l_2 = -(3z_2 X - \ln z_2 - 2)/X,$$
$$m_2 = 2z_2 (\ln z_2 - z_2 X)/X.$$
20. Compute
$$z_3 = z_2 + (1/2)(l_2 + \sqrt{l_1^2 + 4m_2}).$$
21. Compute
$$l_3 = -(3z_3 X - \ln z_3 - 2)/X,$$
$$m_3 = 2z_3 (\ln z_3 - z_3)/X.$$
22. Compute
$$z_4 = z_3 + (1/2)(l_3 + \sqrt{l_n^2 + 4m_n}).$$
23. Compute
$$l_4 = -(3z_4 X - \ln z_4 - 2)/X,$$
$$m_4 = 2z_4 (\ln z_4 - z_4)/X.$$
24. Compute
$$z_5 = z_4 + (1/2)(l_4 + \sqrt{l_4^2 + 4m_4}).$$
25. Compute $\ln(z_5)$. Let it be $Y_1$.
26. Compute
$$z_2' = z_1 + (1/2)(l_1 - \sqrt{l_1^2 + 4m}).$$
27. Compute
$$l_2' = -(3z'_2 X - \ln z'_2 - 2)/X,$$
$$m_2' = 2z_2' (\ln z'_2 - z'_2 X)/X.$$
28. Compute
$$z_3' = z'_2 + (1/2)(l'_2 - \sqrt{l'^2_2 + 4m'_2}).$$
29. Compute
$$l_3' = -(3z'_3 X - \ln z'_3 - 2)/X,$$
$$m_3' = 2z_3' (\ln z'_3 - z_3 X)/X.$$
30. Compute
$$z_4' = z'_3 + (1/2)(l'_3 - \sqrt{l'^2_3 + 4m'_3}).$$

31. Compute
$$l'_4 = -(3z'_4 X - \ln z'_4 - 2)/X,$$
$$m'_4 = 2z'_4 (\ln z'_4 - z'_4 X)/ X.$$
32. Compute
$$z'_5 = z'_4 + (1/2)(l'_4 - \sqrt{l'^2_4 + 4m'_4}.$$
33. Compute $\ln(z_5)$ and $ln(z'_5)$ let these be $Y_1, Y_2$.
34. Print $Y_1, Y_2$ and go to 38.
35. Print $Y = 42.306755092$ and go to **38**.
36. Print 0 and go to 38.
37. Print 'No solution in real domain.'
38. **End**.

This algorithm has been devised with the feature of four iterations for correction, and the formulae used are based on the first method using Equations (3.9), (3.10) and (3.8). Algorithms based on the multiple methods provided in the paper can be easily framed using hints from this algorithm, but they are not explicitly written in the paper for the sake of brevity.

Although this algorithm is meant for the general solution of Lambert W functions, it also addresses a specific case, $W(10^{20}) = y$, with an initial assumption of $z_1 = 1$. This proves its versatility compared to Newton's and Halley's methods, which require initial assumptions within 10% to 20% of the actual values, thus necessitating guesses before solving for a precise approximation. However, this method has solved the function even with an initial assumption that deviates significantly $\{2.363688732(10^{18})\, p.c.\}$ from the correct value.

**6.1. Pseudocode for Solving Lambert W Function**

1. Given x in equation $W(x) = y$.
2. If $x < 0$, go to step 16.
3. If $x = 0$, go to step 38.
4. Initialise $z_1$ to 1.
5. Define the formulae:
$$l_n = -\{3 z_n \ln(z_n) + 2 z_n - x\}/\{\ln(z_n) + 2\},$$
$$m_n = 2 z_n \{x - z_n \ln(z_n)\}/\{\ln(z_n) + 2),$$
$$z_{n+1} = z_n + (1/2)\{l_n + (l_n^2 + 4 m_n)^{1/2}\}.$$
6. Compute $l_1$ and $m_1$ using $z_1$.
7. Compute $z_2$ using $l_1, m_1,$ and $z_1$.
8. Compute $l_2$ and $m_2$ using $z_2$.
9. Compute $z_3$ using $l_2, m_2,$ and $z_2$.
10. Compute $l_3$ and $m_3$ using $z_3$.
11. Compute $z_4$ using $l_3, m_3,$ and $z_3$.
12. Compute $l_4$ and $m_4$ using $z_4$.
13. Compute $z_5$ using $l_4, m_4,$ and $z_4$.
14. Compute $\ln(z_5)$. Let it be Y.
15. Go to step 37.
16. Set X to $-x$. If $1/e < X < \infty$, go to step 39.
17. Initialise $z_1$ to 2.
18. Define the formulae:

$$l_n = -\{3 z_n X - \ln(z_n) - 2\}/X,$$
$$m_n = 2 z_n\{\ln(z_n) - z_n X\}/X,$$
$$z_{n+1} = z_n + (1/2)\{l_n + (l_n^2 + 4 m_n)^{1/2}\}.$$

19. Compute $l_1$ and $m_1$ using $z_1$.
20. Compute $z_2$ using $l_1, m_1$, and $z_1$.
21. Compute $l_2$ and $m_2$ using z2.
22. Compute $z_3$ using $l_2, m_2$, and $z_2$.
23. Compute $l_3$ and $m_3$ using $z_3$.
24. Compute $z_4$ using $l_3, m_3$, and $z_3$.
25. Compute $l_4$ and $m_4$ using $z_4$.
26. Compute $z_5$ using $l_4, m_4$, and $z_4$.
27. Compute $\ln(z_5)$. Let it be $Y_1$.
28. Compute $z_2'$ using $l_1$, $m_1$ and $z_1$:
$$z_2' = z_1 + (1/2)\{l_1 - (l_1^2 + 4 m_1)^{1/2}\}$$
29. Define the formulae:
$$l_n' = -\{3 z_n' X - \ln(z_n') - 2\}/X,$$
$$m_n' = 2 z_n'\{\ln(z_n') - z_n' X\}/X,$$
$$z_{n+1}' = z_n' + (1/2)\{l_n' - (l_n'^2 - 4 m_n')^{1/2}\}.$$

29. Compute $l_2'$ and $m_2'$ using $z_2'$.
30. Compute $z_3'$ using $l_2', m_2'$, and $z_2'$:
$$z_3' = z_2' + (1/2)\{l_2 - (l_2'^2 + 4 m_2')^{1/2}\}$$
31. Compute $l_3'$ and $m_3'$ using $z_3'$.
32. Compute $z_4'$ using $l_3', m_3'$, and $z_3'$:
33. Compute $l_4'$ and $m_4'$ using $z_4'$.
34. Compute $z_5'$ using $l_4', m_4'$, and $z_4'$:
35. Compute $\ln(z_5')$. Let it be $Y_2$.
36. Print y equal to $Y_1$ and $Y_2$, then go to step 40.
37. Print $y$ equal to Y and go to step 40.
38. Print 0 and go to step 40.
39. Print 'No solution in real domain.'
40. End.

## 7. Application of Lambert W Function

Lambert W Function assumes importance because it facilitates solution to many physical processes, where equations are transformable to Lambert W Function. Kindly peruse introduction section. Some equations transformable to Lambert W Function are taken up below to illustrate its applications.

### 7.1. Determination of Of Value of y, When Value of m Is Given in Equation $y^y = m$

Taking natural logarithm of $y^y = m$ and putting $z = \ln y$, $x = \ln(m)$ transforms, $y \ln (y) = \ln (m)$ to $z e^z = x$ or $W(x) = z$, when $x$ lies between 0 and infinity or between 0 and $-1/e$.

### 7.2. Determination Of Value of y, When Value of m Is Given In Equation $y^{1/y} = m$

Substituting $z = \ln y$, $x = \ln(m)$, yields $W(-x) = -z$ which is solvable when $-x$ lies between 0 and infinity or between 0 and $-1/e$.

### 7.3. Determination Of Value of x, When Values of $p, q, r$ Are Given In Equation $p \ln x + q/x = r$

Substitution of $x = y(q/p)$, transforms the equation, $p \ln x + q/x = r$ to $\ln y + 1/y = r/p - \ln(q/p)$, which can be written as $y\, e^{1/y} = (p/q)e^{r/p}$ or $-(1/y)\, e^{-1/y} = -(q/p)e^{-r/p}$ or $W(-X) = -1/y$ where $X = (q/p)e^{-r/p}$ and $-X$ lies between 0 and infinity or between 0 and $-1/e$.

### 7.4. Determination Of Value Of x, When Real Values of $p, q, r$ Are Given In Equation $p \ln x + qx = r$

Substitution of $x = y(p/q)$, transforms the equation $p \ln x + qx = r$ to $\ln y + y = r/p - \ln(p/q)$, which can be written as $y\, e^y = (q/p)e^{r/p}$ or $W(X) = y$ where $X = (q/p)e^{r/p}$ and $X$ lies between 0 and infinity or 0 and $-1/e$.

### 7.5. Determination Of Value of x, When Values of $p, q, r$ Are Given In Equation $p x + q e^{rx} = s$

Substituting $z = e^{rz}$, equation $p x + q e^{rx} = s$, transforms to $(p/r)\ln z + qz = s$ and this equation is same as solved in section 7.4.

### 7.6. Determination Of Value of y, When Value Of x Is Given In Equation $y = x^{x^{x^{x^{\cdot^{\cdot^{\cdot^{x}}}}}}}$

Equation $y = x^{x^{x^{x^{\cdot^{\cdot^{\cdot^{x}}}}}}}$ can be written as, $y = x^y$ or $y^{\frac{1}{y}} = x$. This equation has been solved in section 7.2.

## 8. Comparison of the Present Method with Newton's and Halley's Method of Approximation

To emphasise the salient features and versatility of this method, I make comparison with Newton's and Halley Methods.

### 8.1. Newton's Method of Approximation

For approximation of a root of a function $f(x)$, Newton's Method assigns a value say $x_0$ to be its solution and then improves it to $x_1$ using equation

$$x_1 = x_0 - f(x_0)/f'(x_n),$$

where $f(x_0)$ and $f'(x_0)$ are values of function and its derivative at $x = x_0$. Iterating it $n$ times, the error is appreciably reduced and the equation in general is written

$$x_{n+1} = x_n - f(x_n)/f'(x_n). \tag{8.1}$$

With this background, I highlight some drawbacks with this method.

#### 8.1.1. Initial assumption

Newton's method is highly dependent on the initial assumption for the value to be approximated. If the initial assumption is far from the actual value, *the method may fail to converge or converge to the wrong approximation.* For example, *if the function has multiple roots, a poor choice of the starting point might lead to convergence to a different root than the one intended.* For example, consider the function

$$f(x) = x^3 - 2x + 2.$$

This function has a real root near $x = -1.769$ and other complex roots. If I make initial assumption of $x = 0$, then applying equation (8.1), I obtain, $x_1 = 1, x_2 = 0, x_3 = 1$. Thus, the approximation oscillates between 1 and 0 and does not converge and fails. That proves its high sensitivity to the choice of the initial assumption.

#### 8.1.2. Non-Convergence for Non-Differentiable Functions

Newton's method requires the function to be differentiable [29]. If the derivative at any point is zero or nearly zero, the method may fail or produce excessive large iterations, leading to divergence. For example, $f(x) = x^{1/3}$, the derivative at $x = 0$ is undefined and tends to infinity. Starting with $x = 0$, the method fails to converge.

#### 8.1.3. Slow Convergence for Poor Initial assumptions

Although Newton's method converges quadratically, this fast convergence occurs only when the initial assumption is close to the actual value. If the assumption is far from the actual value, the convergence can be slow. For example, with the function $f(x) = x^2 - 2$, if the initial assumption is $x = 100$, the method may take many iterations to converge to the correct value.

### 8.1.4. Failure at Critical Points of Zero Derivative

If the derivative of the function $f'(x)$ happens to be zero, the method fails. Function $f(x) = x^3 - 3x^2 + 3x - 1$, has derivative $f'(x) = 3x^2 - 6x + 3$, and when $x$ reaches $1$, then $f'(1) = 0$ and the method fails.

### 8.1.5. Requirement for Close Initial Assumption

There is no set single formula for assessing initial assumption close to actual value. In general, percentage variation of 10 to 20 from actual value in initial assumption can achieve for convergence of approximation, and that calls for adoption of Bisection method, Heuristics and physical insight, Graphical methods. The process becomes cumbersome and untidy on account of initial assumption and then final approximation by iterations.

## 8.2. Halley's Method of Approximation

The iterative formula for Halley's method is

$$x_{n+1} = x_n - \frac{2f(x_n)f'(x_n)}{2\{f'(x_n)\}^2 - f(x_n)f''(x_n)}, \quad (8.2)$$

where $f(x_n)$ is the value of the function at $x_n$, $f'(x_n)$ is the first derivative of the function at $x_n$ and $f''(x_n)$ is the second derivative of the function at $x_n$. Although this method is a third-order convergence (cubic) method, meaning it converges faster than Newton's method in theory. However, it has its own set of drawbacks.

### 8.2.1. Complexity and Computational Cost

This method is complicated as it uses second order derivative besides first order derivative of the function as is explicit from the Equation (8.2) and that requires additional second derivative which at times is difficult to compute or expensive to evaluate, especially in higher-dimensional problems. For example, function $f(x) = e^{x^2} sin^2 x$, is a bit difficult and tedious for computing its second derivative.

### 8.2.2. Non-Convergence Near Singularities

Halley's method can fail to converge if the second derivative is zero or approaches zero. For example, the function $f(x) = x^3 - 2x + 2$, has $f'(x) = 3x^2 - 2$ and $f''(x) = 6x$. At $x = 0$, or near zero, $f''(0) = 0$. When $f''(x)$ is small or zero, the denominator in the Halley iteration formula can cause large iterations, leading to non-convergence or divergence of the method.

### 8.2.3. Unnecessary for Simple Functions

For many simple functions, Halley's method provides no real advantage over Newton's method yet introduces additional computational overhead.

### 8.2.4. More Complex to Implement:

Halley's method is inherently more complicated due to the involvement of second-order derivatives. This makes its implementation more difficult, especially for higher-dimensional systems or functions that are not easily differentiable.

### 8.2.5. Potential for Numerical Instability:

In cases, where higher-order derivatives show significant fluctuations, Halley's method can suffer from numerical instability. The method may introduce round-off errors, especially when working with finite precision, and this can lead to divergence. For example, if $f(x) = 1/x^3$, then $f'(x) = -3/x^4$ and $f''(x) = 12/x^5$. If value of $x$ is assumed zero or near zero, then at these points, values of $f'(x)$ and $f''(x)$ goes extremely high which can introduce rounding errors in finite-precision computations. These rounding errors can propagate through the iterations, causing Halley's method to diverge or produce inaccurate results.

## 8.3. Quadratic Equation Approximation

For approximating value of $y$ in Lambert W Function, $W(x) = y$, also written as $ye^y = x$, $e^y$ is considered equal to $z$ and the resultant equation $z\ln(z) = x$, transforms to a quadratic equation $a_1^2 - l_1 a_1 - m_1 \simeq 0$, on substituting $z$ with $z_1 + a_1$, where $z_1$ is initial assumption and $a_1 \ll z_1$. This quadratic equation has coefficients $l_1$ and $m_1$ given by Equations (3.2) and (3.3). Value of $z_{n+1}$ is approximated by recursive relation,

$$z_{n+1} = z_n + a_n = z_n + (1/2)\left(l_n + \sqrt{l_n^2 + 4m_n}\right), \tag{8.3}$$

where $z_{n+1}$ is *(n+1)th* iterative approximation of $z$, $a_n$ is a root of *nth* iterated quadratic equation, $l_n$ and $m_n$ are given by Equations (3.10) and (3.11). This method has salient feature as listed below.

8.3.1. Insensitivity to Initial Assumptions

This method obviates the choice of initial assumption. A value out of a large spectrum depending upon the precision of the calculating device, can be assigned and in each case, the result within four or five iterations come out accurate to nine decimal points.

i. For given values of x equal to $10^{-5}, 0.1, 0.5, 1.0, 100$ and $10^5$, initial assumption of value of $z = 1$ has been made in each case and number of iterations varying from two to four resulted in precise approximation of $y = \ln(z_n)$ up to ten decimal points with zero error. Kindly refer to Table 3.1.

ii. Referring to Table 5.1, for approximation of y in equation $ye^{-y} = x$, when the given value of x was $10^5$ in each case, initial assumptions of y equal to $1, 10, 10^2, 10^3, 10^4, 10^5, 10^6, 10^{12}$, have been made, that is initial assumptions of y has been vastly varied from 1 to $10^{12}$, but still the method led to precise approximation of y without any error with in three to four iterations. That conclusively proves, this method is immune and insensitive to initial assumptions of y.

I could not go beyond $10^{500}$ with initial assumption of $y = 1$ on account of limitation of precision of *Desmos* scientific calculator used for approximation otherwise the limit of $x$ would have gone beyond $10^{500}$ with same initial assumption of $y$ as 1 leading to precise approximation of $y$.

8.3.2. *No Requirement for Differentiability*

A significant feature of this method is that, unlike Newton's or Halley's methods, it does not require the differentiability or continuity of the function. Since no differentiation is involved, it eliminates the time-consuming and complex process of determining first or second-order derivatives of the function.

8.3.3. *Successful at Critical Points of Zero Derivative*

This method does not impose the condition that the derivative of the function at the assigned or iterated value must not be zero. It has an advantage over Newton's and Halley's methods, which fail when the first or second-order derivative of the function is zero or tends towards zero.

8.3.4. *Robustness*

The ability to handle wildly incorrect initial assumptions without failing or slowing down makes this method extremely versatile. This robustness is particularly valuable in real-world applications where initial assumptions are difficult to ascertain, or where functions behave unpredictably.

8.3.5. *Wider Applicability*

As highlighted in section 9.2, this method can be applied to a broader range of problems without requiring preprocessing steps to refine initial guesses—serious drawbacks that affect existing approximation methods.

8.3.6. *Self-Correction and Reliability*

This method has a standalone self-correcting feature, which is particularly valuable as it minimises human error. In contrast, methods such as Newton's or Halley's can be derailed by errors in initial assumptions or calculation setup, leading to wasted time and resources. A method that automatically corrects these errors adds a layer of reliability.

8.3.7. *Fast Convergence*

Converging within 3 to 4 iterations, regardless of the initial guess, is a major advantage over existing methods, many of which require multiple iterations to achieve a high degree of accuracy.

## 9. Results And Conclusions

The crux of the numerical evaluation of Lambert W Function taken up in this paper is linear approximation of $ln(y)$ leading to quadratic approximation of $y \, ln(y)$. The value of $ln(y)$ has been proved to approximate with $2/(2y - 1 - 1/y^3)$ and when value of $y$ is large and $1/y^3$ is negligible, $ln(y)$ approximates with $2/(2y - 1)$ [28]. If $y$ equals, $y_1 + a_1$, then $ln \, y = ln \, y_1 + ln(1 + a_1/y_1)$ approximates $ln \, y_1 + 2/(2y_1/a_1 + 1)$ provided $y_1/a_1$ is large. That needs, assignment of value of $y_1$ in the vicinity of $y$ or $y - y_1 = a_1$ should be exceedingly small. In that case, $y \, ln \, y = x$ approximates with $(y_1 + a_1)ln \, y_1 + 2(y_1 + a_1)/(2y_1/a_1 + 1)$. That is $(y_1 + a_1)ln \, y_1 + 2(y_1 + a_1)/(2y_1/a_1 + 1) - x$ approximates with zero. This is a quadratic equation in $a_1$ which gives two values of $a_1$. Out of these two values, one that yields $y_1 + a_1$ positive quantity is considered. If both roots yield $y_1 + a_1$ positive quantities, then both are considered. When $y - y_1$ is small on account of assignment of value of $y_1$, then approximation of value of $y = y_1 + a_1$ is made in one go.

However, it may not be possible to assign value to $y_1$ so that $y - y_1$ is small, in such cases, $y$ will roughly approximate $y_1 + a_1$. Considering $y = y_1 + a_1 + a_2$, value of $a_2$ is determined by a root of subsequent quadratic equation. Continuing in this fashion, it can be reiterated $n$ times to yield, $y = y_{n+1} = y_n + a_n = y_1 + a_1 + a_2 + a_3 + \cdots + a_n$ or $y = y_1 + \sum_{i=1}^{n} a_i$.

Value of $n$ which is number of iterations, depends upon the precision of result required and it also depends upon the assumption of value of $y_1$. More precise the result required, larger the value of integer $n$ would be. Also, better the proximity of $y_1$ with $y$, fewer iterations would be required. In this paper, in three iterations i.e. = 3, precise results of $y$ could be determined, therefore, further iterations were not required. It was also observed, except in the number of iterations, variation in assigned value to $y_1$ with in certain range did not matter. For example, for given $x = 10^5$ interestingly, variation from 1 to $10^{12}$ in assigned value of $z_1$ did not make any difference in the calculated value of $y$ except the fact when value of $y_1$ is quite away from actual $y$, number of iterations increased.

When $\infty > x \geq 0$, two methods have been used in approximation of $ye^y = x$. In the first method, it is assumed, $z = e^y$ and that transforms the equation to $z \, ln(z) = x$. Considering $z = z_1 + a_1$, where $z_1$ is the initial approximation of z and $a_1 \ll z_1$, the equation $z \, ln(z) = x$ transforms to a quadratic equation in $a_1$ and its roots provide the value of $a_1$ and that results in determining the value of $z_1 + a_1$. For precise approximation, the process is repeated $n$ times, when $z_{n+1} = z_n + a_n$ is achieved and in that situation $z_{n+1}$ converges exponentially to value of $z$.

In second method, $y \, e^y = x$ is written as $ln(y) + y = ln(x)$ and, then this equation transforms to a quadratic equation in $a_1$ considering $y = y_1 + a_1$, where $y_1$ is the initial approximation of $y$ and $a_1 \ll y_1$. Roots of the equation provide the values of $a_1$ and that results in value of $y_1 + a_1$. For precise approximation, the process is repeated $n$ times, when $y_{n+1} = y_n + a_n$ is achieved and in that situation $y_{n+1}$ converges exponentially to exact value of $y$.

For numerical approximation of $y$ in $W(-x) = -y$ or $-ye^{-y} = -x$, when $0 \leq x \leq 1/e$, same two methods as explained in foregoing paragraphs are applicable with a slight difference, now the equation is $-y\,e^{-y} = -x$ and both the roots of the quadratic equations are of positive sign making these useable for approximating two values in stead of one for $y$. In such cases, $y_{n+1}$ and $y'_{n+1}$ both reach precise approximation of $y$. To avoid repetition, the procedure is not reiterated.

### 9.1. Self-Corrective Nature Of Quadratic Approximation

Because of wide range of values that can be assigned to $y_1$, inadvertent error in calculations is self-corrected by subsequent iterations provided the error is not fatal throwing the value of $y_1$ or $z_1$ as the case may be, out of permissible range. Inherent self correction present in this type of quadratic approximation is akin to inherent presence of negative feedback in stable electronic system. Needless to say, self corrective nature of approximation, in no way, calls for lackadaisical approach on the part of the operator carrying out computations.

It is submitted that negative or complex values of $y_1 + a_1$ or $z_1 + a_1$ as the case may be, were ignored. The question that arises is what had happened if these were duly considered in our evaluation. Interestingly, had these been considered, these would have resulted in complex values of $y$. It was also considered, in equation $-ye^{-y} = -x$, value of $x$ should be such that $0 < x \le 1/e$. Had the value of $x$ been more than $1/e$ in equation, $-ye^{-y} = -x$, the results in both the cases would have been complex and multiple in numbers. This is the reason Lambert W Function has multiple solutions. Exhaustive examples provided in tabular form in the paper, leaving nothing to imagination have proved the veracity of all the formulae derived. The paper dealt with empirical and verifiable research by providing unique method of quadratic equation approximation. This method is different from *Newton method of approximation* or *Halley method of approximation,* and it provided precise results in as few as three iterations. It also provided algorithm for numerical evaluation of $y$ in Lambert W Function $W(x) = y$, where $\infty > x \ge -1/e$.

### 9.2. Potential of Quadratic Equation Approximation for Future Use in Research

It is well known that many processes in the physical and life sciences make use of transcendental equations involving exponential and logarithmic terms. Additionally, some processes rely on double exponentials. For example, in theoretical biology for population growth, number theory for Fermat numbers, double Mersenne numbers, Sylvester's sequence, and *k-ary* Boolean functions. In computer science, double exponentials are employed in algorithmic complexities, and in physics, for modelling systems like the Toda oscillator or dendritic macromolecules. Similarly, logarithmic, and double logarithmic functions (i.e., the logarithm of a logarithm) are also integral to many physical processes.

Quadratic equation approximations can be applied to estimate the variables and outputs of processes involving such functions. This approach is particularly relevant when using inverse functions of the form $N(X) = y$, where $ye^{e^y} = X$ as applied to various physical phenomena. A research paper discussing this has been written, accepted, and is soon to be published. Notably, this inverse function can not be directly solved using the Lambert W function, and there appears to be no straightforward solution to approximate its value for a given real input. The method proposed in this paper has also been successfully applied to the solution of the Lambert W function $W(a + ib) = x + iy$, in the complex domain, where $(x + iy)e^{x+iy} = a + ib$, and $a$ and $b$ are given real quantities, This method of approximation shows potential for solving polynomial equations, particularly
$$x^{2n+1} - a\,x = b,$$
when $a$ has a positive real value, $b$ is real, and $n$ is a positive integer.

### 10. Author's Contribution

Under this head, I would like to emphasise the following key aspects of my methodology.

The approximation of $\ln\{(y)/(y-1)\}$ by $2/(2y-1)$, when $y$ is large, was developed and published in an earlier article 'Approximation Of Logarithm, Factorial And Euler - Mascheroni Constant Using Odd Harmonic Series, in journal *Mathematical Forum*. This paper is also listed at serial number 23 in the references. Since the Lambert W function is the inverse of the equation $ye^y = x$, where $x$ is a given real quantity, this relationship can be written as $y + ln(y) = ln(x)$. By substituting $y = y_1 + a_1$, where $y_1$ is an initial approximation, and applying my approximation of $ln(y)/(y - 1)$, the equation was transformed into a quadratic equation in terms of $a_1$.

The roots of this quadratic equation were then used to approximate the value of $y = y_1 + a_1$. By iterating this process using the formula $y_{n+1} = y_n + a_n$, the solution converged to a highly precise value for $y$. This iterative approach, rooted in a quadratic equation, was named *'Quadratic Equation Approximation.'*

The concept and methodology of this quadratic approximation are entirely original and have been developed solely by me. I have contributed 100% to this paper, including the derivation of the approximation technique and the iterative refinement process.

In addition to the work presented in this paper, the numerical approximation using quadratic equation approximation of the special function $e^{e^y} = X$, where $X$ is a given real value, has also been accepted for publication in another journal. It is important to note that this special function is not solvable by the Lambert W Function, further emphasizing the novelty of the proposed approximation method. This method also has the potential of approximating a real root of odd degree polynomial equations.

**Statements and Declarations**

*Competing Interests, Conflict of Interest, Funding and* **Acknowledgement***:* There are no competing interests financial or non-financial involved in this research work. Further, there is no conflict of interest of any kind involved in this research work. No fund whatsoever in cash or kind has been received for this work. However, I acknowledge the help provided by the website https://www.desmos.com in calculating the values of tedious exponential and logarithmic terms.

**Data Availability Statement**
There is no data availability statement.

**Ethical Approval**
Not applicable.